\documentclass[aop,MSNbibl,seceqn,citesort,dvips]{arximspdf}
\usepackage{mathrsfs}

\doi{10.1214/11-AOP667}
\volume{40}
\issue{5}
\pubyear{2012}
\firstpage{2131}
\lastpage{2167}

\makeatletter
\newtheorem{Th}{Theorem}[section]
\newtheorem{lm}[Th]{Lemma}
\newproclaim{re}[Th]{Remark}

\setattribute{abstract}{width}{285pt}
\setattribute{abstract}{skip}{8\p@}

\makeatother

\begin{document}
\begin{frontmatter}

\title{Diffusion processes in thin tubes and~their~limits~on~graphs}
\runtitle{Diffusion processes in thin tubes}

\begin{aug}
\author[A]{\fnms{Sergio} \snm{Albeverio}\thanksref{t1}\ead[label=e1]{albeverio@uni-bonn.de}}
\and
\author[B]{\fnms{Seiichiro} \snm{Kusuoka}\corref{}\thanksref{t2}\ead[label=e2]{kusuoka@math.kyoto-u.ac.jp}}
\runauthor{S. Albeverio and S. Kusuoka}
\affiliation{Universit{\"a}t Bonn and Kyoto University}
\address[A]{Institut f\"ur Angewandte Mathematik\\
Universit{\"a}t Bonn\\
Endenicherallee 60\\
 53115, Bonn\\
Germany\\
\printead{e1}} 
\address[B]{Graduate School of Science\\
Kyoto University\\
Kitashirakawa-Oiwakecho\\
Kyoto 606-8264\\
Japan\\
\printead{e2}}
\end{aug}
\thankstext{t1}{Supported by Provincia Autonoma di Trento, through the
NEST-Project.}
\thankstext{t2}{Supported by the Excellent Young Researcher Overseas Visit Program.}

\received{\smonth{10} \syear{2010}}
\revised{\smonth{2} \syear{2011}}

%
\begin{abstract}
The present paper is concerned with diffusion processes running on
tubular domains with conditions on nonreaching the boundary,
respectively, reflecting at the boundary, and corresponding processes in
the limit where the thin tubular domains are shrinking to graphs.
The methods we use are probabilistic ones.
For shrinking, we use big potentials, respectively, reflection on the
boundary of tubes.
We show that there exists a unique limit process, and we characterize
the limit process by a second-order differential generator acting on
functions defined on the limit graph, with Kirchhoff boundary
conditions at the vertices.
\end{abstract}

%

\end{frontmatter}

\section{Introduction}

The present paper is concerned with diffusion processes running on
tubular domains with Dirichlet (i.e., absorbing-like) (resp., Neumann,
i.e., reflecting) boundary conditions, and the respective processes
obtained in the limit where the thin tubular domains shrink to graphs.
Problems of this type have been intensively studied before in the case
of Neumann boundary conditions, both by probabilistic tools~\cite{FW,FW3} and analytic tools
~\cite{ACF,Bon,CE,CF,DC,DT,EP,MV,Po}.
The case of Dirichlet boundary conditions was known to present special
difficulties, which explains why there have been, up to now, fewer
works concerned with this case, and, in fact, these are only concerned
with either special graphs or special shrinking procedures, leading
mainly (with the exception of~\cite{ACF,CE,CF,DC}) to limiting processes which ``decouple at
vertices''
\cite{EP,BM,CMN}.

Before explaining these difficulties and entering into details let us
motivate the reasons to undertake such studies, pointing out also some
connections with other problems and giving some historical
remarks.\vadjust{\goodbreak}

In many problems of analysis and probability one encounters
differential operators defined on structures which have small
dimensions in one or more directions.
Let us mention as examples the modeling of fluid motion in narrow
tubes, or in nearly two-dimensional domains (see, e.g.,~\cite{Ra}), the
propagation of electric signals along nearly one-dimensional neurons
(see, e.g.,~\cite{ADM,BM,CMN}), the propagation of
electromagnetic waves in wave guides~\cite{KMS}, the propagation of
quantum mechanical effects in thin wires (in the context of
nanotechnology); see, for example,~\cite{ACF,CE,CF,DC,DT,EP,ES,GP,Kuc1,Kuc2,Kur,Po,Tu}.
Such geometrical structures tend in a certain limit (mathematically
well described in general through a Gromov topology) to a graph.
Modeling dynamical systems or processes on such structures by
corresponding ones on a graph might present certain advantages (e.g.,
PDEs becoming ODEs on graphs; more dimensional spectral problems
reduced to one-dimensional ones). In any case the study of dynamics and
processes on graphs can be considered as an idealization or a ``first
approximation'' for the study of the corresponding objects in more
realistic situations.

There is a rich literature on differential operators on graphs.
Diffusion operators and evolution equations were considered originally
in work by Lumer~\cite{Lu}, and subsequently by many authors; see, for
example,~\cite{vB,Ya2,AvN,Mu2}.
Elliptic and parabolic nonlinear equations on graphs have been
discussed, for example, in relations to applications in biology, for
example, in~\cite{CMN}; see also, for example,~\cite{BM,ADM}
for nonlinear diffusions on graphs in connection with neurobiology.
Heat kernels on graphs have been studied in particular in~\cite{Mu}.
Hyperbolic nonlinear equations on graphs have been studied, for
example, in~\cite{KMS}.

In quantum mechanics, Schr\"odinger equations on graphs are considered
as models of nanostructures; see, for example,~\cite{ES,BCFK,Kuc1,Kuc2}.
Work has been particularly intense in the study of spectral properties
of Sch\"odinger-type operators on graphs; see, for example,~\cite{GP,Kuc1,Kuc2,Kur}.
Such models of quantum mechanics on graphs also play an important role
in the study of the relation between classical chaos and quantum chaos;
see, for example,~\cite{Kur,GP,ES2,Sa,RS}.

For the study of the limit of differential operators on thin domains of
${\mathbb R}^n$ (and corresponding PDEs) degenerating into geometric
graphs (and corresponding ODEs) we refer to~\cite{Ya2,Ra,Ko} and especially to the surveys by Raugel~\cite{Ra} (which discuss
topics like spectral properties, asymptotics and attractors).
For the study of parabolic equations and associated semi-groups and
diffusion processes we also refer to~\cite{Ra}.
Corresponding hyperbolic problems in connection with the modeling of
ferroelectric materials have been discussed, for example, in~\cite{AH}.

Probabilistic methods for the study of processes on thin domains of
${\mathbb R}^n$ have been developed by Freidlin and Wentzell in the
case of Neumann boundary conditions.
They exploit the consideration of slow, respectively fast, components
going back to~\cite{FW2}, applied to the thin tubes problem~\cite{FW}.
In these studies the basic probabilistic observation is that for a
Brownian motion in a thin tube along a line, the component in the
transverse direction is fast, and the one in the longitudinal direction
is slow.
The control in the limit exploits the assumption on the reflecting
properties of the fast component, together with a projection technique
onto the longitudinal direction.
In~\cite{FW} it is shown that the diffusion coefficient for this limit
process is obtained by averaging the diffusion coefficient for the
process in tubular domains with respect to the invariant measure of the
fast component with suitable changed space and time scales.

Analytically the Laplacian in the transverse direction has a constant
eigenvalue $0$ (ground state in the transverse direction), which then
yields a natural identification of the subspace of $L^2$---over the
thin tube corresponding to the eigenvalue $0$ for the Laplacian in the
transverse direction with the $L^2$---space along an edge.
Results about this approximation concern convergence of eigenvalues,
eigenfunctions, resolvents and semigroups~\cite{EP,Gr,MV,DT}.
Besides, operatorial and variational methods also methods of Dirichlet
form theory have been used~\cite{Bon}.

The identification stressed above is no longer possible in the case of
Dirichlet boundary conditions on the boundary of the thin tube, since
the lowest eigenvalue of the Laplacian in the transverse direction
diverges like $1/\varepsilon^2$, where $\varepsilon>0$ is the width
of the narrow tube.
(For a probabilistic study of the first-order asymptotics of the lowest
eigenvalue of the Dirichlet Laplacian in tubular neighborhoods of
submanifolds of Riemannian manifolds, see~\cite{KP}.)
This has been pointed out clearly and posed as an open problem by Exner
(see~\cite{AGHH}).
In order to nevertheless manage analytically the limit to a graph, one
has to perform a renormalization procedure, first introduced in \cite
{ACF}, and extended in~\cite{CE,CF}, for the case of a V-graph
(waveguide).
More general cases with Dirichlet boundary conditions have been managed
in the case where the shrinking at vertices is quicker than the one at
the edges; however, then one has ``no communication between the
different edges'' (i.e., ``decoupling'') on the graphs; see~\cite{Gr,MV,Po}.
The interest in discussing the case of Dirichlet-boundary conditions is
particularly clear in the physics of conductors, where such boundary
conditions arise most naturally, both in classical and quantum
mechanical problems.
However, in the other type of applications we have mentioned there is
also an interest in studying boundary conditions that are different
from the Neumann ones, since boundary conditions influence the limit
behavior, and one is interested to obtain on the graphs the most
general possible boundary conditions at the vertices (even in the case
of an ``$N$-spider graph'' there are $N^2$-different possible
self-adjoint realizations of a Laplacian on the spider; see, for
example,~\cite{ES,KS2}).

The present paper mainly discusses the case of shrinking by potentials,
and the goal is to determine the limit process on a given graph.
This shrinking by potentials corresponds to confining the process in
thin tubes around the graph, not reaching the boundary almost surely,
and in this sense is related with Dirichlet boundary conditions (the
latter property corresponding however to a completely absorbing boundary).
In Sections~\ref{section curve} and~\ref{section spider} we consider
special cases, because the consideration of these cases illustrate
better the methods we use.

In Section~\ref{section curve} the case of a thin tube $\Omega
^\varepsilon$ in ${\mathbb R}^n$ shrinking to a curve $\gamma$ in
${\mathbb R}^n$ is discussed.
The tube $\Omega^\varepsilon$ has a uniform width $\varepsilon>0$.
In the tube we have a~nondegenerate diffusion process $X^\varepsilon$
with a drift consisting of two parts, one continuous and bounded, the
other of gradient type, pushing away from the boundary, so that the
first hitting time of $X^\varepsilon$ at the boundary $\partial\Omega
^\varepsilon$ is infinite almost surely.
We also construct a diffusion process $X$ on $\gamma$ and show (Theorem
\ref{th-curve}) that if $X^\varepsilon(0)$ converges weakly to $X(0)$,
then also $X^\varepsilon$ converges weakly to $X$.\
If pathwise uniqueness holds both for $X^\varepsilon$ and $X$,
then~$X^\varepsilon$ also converges to $X$ almost surely as $\varepsilon
\downarrow0$.
We also state corresponding results for a process in $\Omega
^\varepsilon$ with a reflecting boundary condition on the
boundary~$\partial\Omega^\varepsilon$ (Theorem~\ref{th-curve2}).
These results are obtained in a similar way as those obtained by our
shrinking with potentials in the first part of Section~\ref{section
curve}.\looseness=-1

In Section~\ref{section spider} we discuss the case of shrinking $N$
thin tubes in ${\mathbb R}^n$ to an $N$-spider graph in ${\mathbb R}^n$.
In this section, we often use the methods discovered by Freidlin and
Wentzell~\cite{FW}, extend their method to the case of diffusion
processes instead of Brownian motions and apply it to the case of
shrinking by potentials.
The process $X^\varepsilon$ in the domain $\Omega^\varepsilon$
consisting of $N$ tubes is defined in a similar way as in Section \ref
{section curve}, $\varepsilon>0$ being the parameter of shrinking to
the $N$-spider graph $\Gamma$ for $\varepsilon\downarrow0$.
We prove again that the first hitting time of $X^\varepsilon$ at the
boundary $\partial\Omega^\varepsilon$ is infinite and that the laws
of $\{ X^\varepsilon\dvtx \varepsilon>0\}$ are tight in the topology of
probability measures on $C([0,+\infty))$, if their initial
distributions are tight.
We then show that any limit process is strong Markov and study the
transition probabilities from the vertex $O$ to any edge of the spider
graph $\Gamma$.
This requires quite detailed estimates of the behavior of the process
$X^\varepsilon$ in a neighborhood of $O$ in $\Omega^\varepsilon$.
These results imply that the boundary condition at $O$ should be a
weighted Kirchhoff boundary condition for the functions in the domain
of the generator of the limit processes $X$.
(This is one of the types of boundary conditions known from the general
discussions on boundary conditions for processes on graphs; see, for
example,~\cite{KS2,Kuc1,Kuc2,KZ,ES,DC}.)
The weights are determined explicitly from the construction, as
transition probabilities to the edges (Lemma~\ref{lem3.7}).
This is crucial to determine the generator of the unique limit process
$X$ (Theorem~\ref{th-spider}).
Similar considerations lead to corresponding results for the case where
$X^\varepsilon$ is a diffusion in $\Omega^\varepsilon$ with reflecting
boundary conditions on $\partial\Omega^\varepsilon$ (Theorem \ref
{th-spider2}).

In Section~\ref{section graph} we state the results in the case of thin
tubes around general graphs, which are obtained immediately from the
results in Sections~\ref{section curve} and~\ref{section spider}.
These are systems consisting of thin tubes around finitely ramified
graphs in ${\mathbb R}^n$ with edges which consist of $C^3$-curves.
Theorem~\ref{th-graph} presents a result similar to the one for an
$N$-spider graph, showing, in particular, convergence of the diffusion
process $X^\varepsilon$ not leaving the system $\Omega^\varepsilon$ of
tubes around the general graph to a diffusion process $X$ on the graph.
Again its generator is determined and an extension is given to the case
of a diffusion with reflecting boundary conditions on $\partial\Omega
^\varepsilon$.
Since the latter result is not only for a~Brownian motion in the thin
tubes, but also for reflecting diffusion processes in the thin tubes,
it is also an extension of previous results of Freidlin and Wentzell
\cite{FW}.

All random variables discussed in the present paper are defined on a
probability space with probability measure $P$, and $E[\cdot]$ denotes
their expectation with respect to $P$.
For a locally compact topological subspace $A$ of ${\mathbb R}^n$, let
$C_0(A) := \{ f\in C(A)\dvtx  \lim_{|x|\rightarrow+\infty} f(x)=0 \}$.

\section{The case of curves}\label{section curve}

In this section, we consider shrinking of thin tubes to curves.
Let $n$ be an integer larger than or equal to $2$.
Let $\gamma\in C^3 ({\mathbb R};{\mathbb R}^n)$ such that $ |\dot
\gamma|=1$ [with $\dot\gamma$ the derivatives of $t \rightarrow
\gamma(t)$, and $|\cdot|$ the norm in ${\mathbb R}^n$], and assume
that $\gamma$ has no self-crossing point, and $\ddot\gamma$ is
a bounded function with a compact support. Let $\varepsilon>0$, $\langle\cdot,\cdot\rangle$ be the
inner product on ${\mathbb R}^n$, and $d(x ,\gamma)$ be the distance
between $x$ and $\gamma$.
Note that $d(x ,\gamma)$ is Lipschitz continuous in $x$.
Define domains $\{ \Omega^\varepsilon\}$ by
\[
\Omega^\varepsilon:= \{ x\in{\mathbb R}^n\dvtx  d(x,\gamma)< \varepsilon
\} .
\]
Consider a differentiable function $u$ on $[0,1)$ such that
\[
u(0)=0,\qquad u' \geq0 \qquad  \lim_{R \uparrow 1} u'(R) = +\infty \quad\mbox{and}\quad-\lim_{R \uparrow1}
\frac{u(R)}{\log(1-R)} = +\infty.
\]
For example, if we define $u(r):= r^\alpha/ (1-r^\alpha)$ for $r\in
[0,1)$ where $\alpha>0$,
then~$u$ satisfies the conditions above.
Let
\[
U^\varepsilon(x) = u ( \varepsilon^{-1} d(x,\gamma)
),\qquad x \in\Omega^\varepsilon.
\]
For $\varepsilon>0$, consider a diffusion process $X^\varepsilon$
given by the following equation:
%
\begin{eqnarray}\label{SDE1}
X^\varepsilon(t)&=&X^\varepsilon(0) + \int_0 ^{t\wedge\zeta
^\varepsilon}\sigma(X^\varepsilon(s))\,dW(s) + \int_0^{t\wedge\zeta
^\varepsilon}b(X^\varepsilon(s))\,ds \nonumber
\\[-8pt]
\\[-8pt]
&& {}- \int_0^{t\wedge\zeta^\varepsilon}(\nabla U^\varepsilon
)(X^\varepsilon(s))\,ds,
\nonumber
\end{eqnarray}
where $X^\varepsilon(0)$ is an $\Omega^\varepsilon$-valued random
variable, $W$ is an $n$-dimensional Wiener process, $ \sigma\in
C_b({\mathbb R}^n; {\mathbb R}^n \otimes{\mathbb R}^n) $, $ b \in
C_b({\mathbb R}^n; {\mathbb R}^n)$ and $\zeta^\varepsilon$ is the
first hitting time of $X^\varepsilon$ at the boundary $\partial\Omega
^\varepsilon$ of $\Omega^\varepsilon$.
Let $a := \sigma\sigma^T$ (with $\sigma^T$ the transpose of~$\sigma
$), and assume that $a$ is a uniformly positive definite matrix.
Then, the solution $X^\varepsilon$ of~(\ref{SDE1}) exists uniquely;
see, for example,~\cite{SV}.

\begin{lm}\label{lm1-1}
$\zeta^\varepsilon= +\infty$ almost surely for small
$\varepsilon>0$.\vadjust{\goodbreak}
\end{lm}

\begin{pf}
Assume $n\geq3$.
Note that $X^\varepsilon$ does not hit $\gamma$ almost surely in this case.
Let $X_x^\varepsilon$ be the solution of~(\ref{SDE1}) replacing
$X^\varepsilon(0)$ and $\zeta^\varepsilon$ by $x$ and $\zeta
_x^\varepsilon$, respectively, where $\zeta_x^\varepsilon$ is the first
hitting time of $X^\varepsilon_{x}$ at $\partial\Omega^\varepsilon$.
It is sufficient to show that $\zeta_x^\varepsilon=+\infty$ almost
surely for $x$ near to $\partial\Omega^\varepsilon$.
By the tubular neighborhood theorem and Theorem 1 in~\cite{Fo}, there
exists a $C^2$-diffeomorphism $\phi= (\phi_1, \phi_2)$ from $\Omega
^\varepsilon\setminus\gamma$ to $\{ y=(y_1,y_2)\in{\mathbb R}\times
{\mathbb R}^{n-1} \dvtx  0<|y_2|<\varepsilon\}$ which satisfies, for small
$\varepsilon$,
\[
\phi_1 (x) = \gamma^{-1}\circ\pi(x)  \quad \mbox{and} \quad  \phi_2 (x)=
d(x,\gamma)\nabla d(x,\gamma),\qquad x\in\Omega^\varepsilon
\setminus\gamma,
\]
where $\pi(x)$ is the nearest point in $\gamma$ from $x$.
Note that $\phi$ is a $C^2$-function on $\Omega^\varepsilon$ and
$\langle\nabla\pi,\nabla U^\varepsilon\rangle=0$ for small
$\varepsilon$. Hence, $\langle\nabla\phi_1 ,\nabla U^\varepsilon
\rangle=0$ and $\nabla\phi_2 \nabla U^\varepsilon= \varepsilon
^{-1}u'(\varepsilon^{-1}d(\cdot,\gamma))\nabla d(\cdot,\gamma)$. By
It\^o's formula, we have
%
\begin{eqnarray}\label{SDE51}
\phi_1(X_x^\varepsilon(t))&=&\phi_1(x) +
\int_0 ^{t\wedge\zeta_x^\varepsilon} \nabla\phi_1 (X_x^\varepsilon
(s)) \sigma(X_x^\varepsilon(s))\,dW(s)\nonumber \\
&&{}+ \int_0^{t\wedge\zeta_x^\varepsilon} \nabla\phi
_1 (X_x^\varepsilon(s)) b(X_x^\varepsilon(s))\,ds \\
&&{}+ \frac12 \sum_{i,j=1}^n \int_0 ^{t\wedge\zeta
_x^\varepsilon} a_{ij} (X_x^\varepsilon(s))\, \partial_i \partial_j
\phi_1 (X_x^\varepsilon(s))\,ds,\nonumber
\\[-1pt]
\label{SDE52}
\phi_2(X_x^\varepsilon(t))&=&\phi_2(x)+
\int_0 ^{t\wedge\zeta_x^\varepsilon} \nabla\phi_2 (X_x^\varepsilon
(s)) \sigma(X^\varepsilon_x(s))\,dW(s) \nonumber\\
&&{}+ \int_0^{t\wedge\zeta_x^\varepsilon} \nabla\phi
_2 (X_x^\varepsilon(s)) b(X^\varepsilon_x(s))\,ds \nonumber
\\[-8.25pt]
\\[-8.25pt]
&&{}+ \frac12 \sum_{i,j=1}^n \int_0 ^{t\wedge\zeta
_x^\varepsilon} a_{ij}(X_x^\varepsilon(s))\, \partial_i \partial_j
\phi_2 (X_x^\varepsilon(s))\,ds\nonumber \\[-1pt]
&&{}- \varepsilon^{-1}\int_0^{t\wedge\zeta
_x^\varepsilon} u'(\varepsilon^{-1}d(X_x^\varepsilon(s) ,\gamma
))\nabla d(\cdot,\gamma)|_{X_x^\varepsilon(s)}\,ds.\nonumber
\end{eqnarray}
Moreover, again by It\^o's formula,
\begin{eqnarray*}
&& |\phi_2(X_x^\varepsilon(t))|^2 \\[-2pt]
&& \qquad  =|\phi_2(x)|^2 + 2\int_0 ^{t\wedge\zeta^\varepsilon_x
} \langle\phi_2 (X_x^\varepsilon(s)), \nabla\phi_2
(X_x^\varepsilon(s)) \sigma(X_x^\varepsilon(s))\,dW(s) \rangle\\[-0.5pt]
&& \qquad  \quad {} + 2\int_0^{t\wedge\zeta_x^\varepsilon}
\langle\phi_2 (X_x^\varepsilon(s)), \nabla\phi_2 (X_x^\varepsilon
(s)) b(X_x^\varepsilon(s))\rangle \,ds \\[-0.5pt]
&& \qquad  \quad {} + \int_0 ^{t\wedge\zeta_x^\varepsilon}
\Biggl\langle\phi_2 (X_x^\varepsilon(s)), \sum_{i,j=1}^n
a_{ij}(X_x^\varepsilon(s))\, \partial_i \partial_j \phi_2
(X_x^\varepsilon(s))  \Biggr\rangle \,ds \\[-0.5pt]
&& \qquad  \quad {} - 2\varepsilon^{-1}\int_0^{t\wedge\zeta
_x^\varepsilon} |\phi_2 (X_x^\varepsilon(s))|u'(\varepsilon
^{-1}d(X_x^\varepsilon(s) ,\gamma))\,ds\\
&& \qquad  \quad {} +\int_0 ^{t\wedge\zeta_x^\varepsilon} \operatorname{trace}  [ \nabla\phi_2 (X_x^\varepsilon(s)) \sigma
(X_x^\varepsilon(s)) ( \nabla\phi_2 (X_x^\varepsilon(s)) \sigma
(X_x^\varepsilon(s)) ) ^T  ]\,ds.
\end{eqnarray*}
Let
\begin{eqnarray*}
\bar a &:=& \sup \bigl\{  | ( \nabla\phi_2 (x) \sigma(x)) ^T \xi
 | ^2 \dvtx  x\in\Omega^\varepsilon, \xi\in\{ y\in{\mathbb R}^n\dvtx
|y|=1\}  \bigr\},\\
\bar b &:=& \sup_{x\in\Omega^\varepsilon} \Biggl( 2\langle\phi_2
(x), \nabla\phi_2 (x) b(x)\rangle+  \Biggl\langle\phi_2 (x) ,\sum
_{i,j=1}^n a_{ij}(x)\, \partial_i \partial_j \phi_2 (x) \Biggr \rangle
  \\[-3pt]
&&\hphantom{\inf_{x\in\Omega^\varepsilon} \Biggl(}\hspace*{82.7pt} {}+ \operatorname{trace}  [ \nabla
\phi_2 (x) \sigma(x) ( \nabla\phi_2 (x) \sigma(x) ) ^T
 ]  \Biggr) .
\end{eqnarray*}
Take $c_0 \in(0,1)$ such that $\sup _{x\in [c_0,1)} (\bar b - 2xu'(x)) \leq 0$ and
\[
f(x):= \int_{c_0^2\varepsilon^2} ^x \exp \biggl( -2\int
_{c_0^2\varepsilon^2}^y \frac{\bar b -2\varepsilon^{-1}\sqrt z
u'(\varepsilon^{-1}\sqrt z)}{{\bar a} z}\,dz  \biggr)\,dy, \qquad x \in
[0,\varepsilon^2).
\]
Then, by It\^o's formula, for $\delta$ such that $0< \delta< 1-c_0$
and for $x$ such that $ c_0 \varepsilon\leq d(x,\gamma)\leq
\varepsilon(1-\delta)$, we have that
\[
E\bigl[f\bigl(\bigl|\phi_2\bigl(X_x^\varepsilon\bigl(T^{c_0 \varepsilon}\wedge T^{\varepsilon
(1-\delta)}\bigr)\bigr)\bigr|^2\bigr)\bigr] \leq f(d(x,\gamma)^2),
\]
where $ T^c:=\inf\{ t>0\dvtx  d(X_x^\varepsilon,\gamma)= c\}$ for $c>0$.
Since
\begin{eqnarray*}
&& E\bigl[f\bigl(\bigl|\phi_2\bigl(X_x^\varepsilon\bigl(T^{c_0 \varepsilon}\wedge
T^{\varepsilon(1-\delta)}\bigr)\bigr)\bigr|^2\bigr)\bigr] \\
&& \qquad  = f ( c_0^2\varepsilon^2 ) P \bigl( T^{c_0 \varepsilon} <
T^{\varepsilon(1-\delta)} \bigr) + f \bigl( \varepsilon^2 (
1-\delta ) ^2 \bigr) P \bigl( T^{c_0 \varepsilon} > T^{\varepsilon
(1-\delta)} \bigr)
\end{eqnarray*}
and
\[
P \bigl( T^{c_0 \varepsilon} < T^{\varepsilon(1-\delta)} \bigr) +
P\bigl ( T^{c_0 \varepsilon} > T^{\varepsilon(1-\delta)} \bigr) =1,
\]
we have
\[
P \bigl( T^{c_0 \varepsilon} > T^{\varepsilon(1-\delta)} \bigr) \leq
\frac{f(d(x,\gamma)^2) - f ( c_0^2\varepsilon^2 )}{f (
\varepsilon^2(1-\delta) ^2 ) - f ( c_0^2\varepsilon^2 )}.
\]
The assumptions on $u$ imply that $f( \varepsilon^2( 1-\delta) ^2)$
diverges to $+\infty$ as $\delta\rightarrow0$. Hence, the proof is
achieved from the fact that $T^{\varepsilon(1-\delta)}$ converges to
$\zeta_x^\varepsilon$ as $\delta\rightarrow0$.

In the case where $n=2$, since $X^\varepsilon$ can hit $\gamma$, we
need a little arrangement.
Let $\Omega_+ ^\varepsilon$ and $\Omega_- ^\varepsilon$ be the two
domains consisting of $\Omega^\varepsilon\setminus\gamma$, and
$\theta^\varepsilon(x)$ be $1$ if $x \in\Omega_+ ^\varepsilon$, $-1$
if $x \in\Omega_- ^\varepsilon$ and $0$ if $x \in\gamma$.
By the tubular neighborhood theorem and Theorem~1 in~\cite{Fo} again,
there exists a $C^2$-diffeomorphism $\phi= (\phi_1, \phi_2)$ from
$\Omega^\varepsilon$ to $\{ y=(y_1,y_2)\in{\mathbb R}\times
(-\varepsilon, \varepsilon) \}$ which satisfies, for small
$\varepsilon$,
\[
\phi_1 (x) = \gamma^{-1}\circ\pi(x)  \quad \mbox{and} \quad  \phi_2 (x)=
\theta^\varepsilon(x)d(x,\gamma),\qquad x\in\Omega^\varepsilon,
\]
such that~(\ref{SDE51}) and~(\ref{SDE52}) hold.
Thus, we can discuss this case in a similar way as the case where
$n\geq3$.
\end{pf}

\begin{Th}\label{th-curve}
Define a diffusion process $X$ by the solution of the following equation:
%
\begin{eqnarray}\label{SDE10}
  X(t)
    &=&    X(0) + \int_0 ^t \dot\gamma\circ\gamma^{-1}
(X(s)) \langle\dot\gamma\circ\gamma^{-1} (X(s)), \sigma
(X(s))\,dW(s) \rangle \nonumber\\
&&{} + \int_0^t \dot\gamma\circ\gamma^{-1} (X(s)) \langle\dot
\gamma\circ\gamma^{-1} (X(s)), b(X(s)) \rangle \,ds \nonumber\\
&&{} + \frac 12 \int _0 ^t \ddot \gamma \circ \gamma ^{-1} (X(s)) | \sigma (X(s))^T \dot \gamma \circ \gamma ^{-1} (X(s))| ^2 \,ds \\
&&{} + \int _0 ^t \dot \gamma \circ \gamma ^{-1} (X(s))
\langle \sigma (X(s))^T \ddot \gamma \circ \gamma ^{-1} (X(s)),\nonumber\\
&&\hspace*{79pt}\qquad{}\sigma (X(s))^T \dot \gamma \circ \gamma ^{-1} (X(s))\rangle\, ds.\nonumber
\end{eqnarray}
Note that $X$ is uniquely determined as a process on $\gamma$.

If $X^\varepsilon(0)$ converges to a $\gamma$-valued random variable
$X(0)$ weakly, then the process $X^\varepsilon$ converges weakly to $X$
in the sense of their laws on $C([0,+\infty); {\mathbb R}^n)$ as
$\varepsilon\downarrow0$.

Moreover, if pathwise uniqueness holds for~(\ref{SDE10}) and (\ref
{SDE1}) for all $\varepsilon>0$, and $X^\varepsilon(0)$ converges to a
$\gamma$-valued random variable $X(0)$ almost surely,
then~$X^\varepsilon$ converges to $X$ almost surely, as $\varepsilon
\downarrow0$.
\end{Th}

\begin{pf}
Note that equation~(\ref{SDE51}) holds even if we replace
$X_x^\varepsilon$, $x$ and $\zeta_x^\varepsilon$ by $X^\varepsilon$,
$X^\varepsilon(0)$ and $\zeta^\varepsilon$, respectively.
Lemma~\ref{lm1-1} implies
%
\begin{equation}\label{th-curve1-1}
\sup_{t\in[0,+\infty)} d(X^\varepsilon(t), \gamma) \rightarrow
0,\qquad\varepsilon\downarrow0,
\end{equation}
almost surely.
Hence, the boundedness of the coefficients implies the tightness of the
process $\phi_1(X^\varepsilon)$.
Let $X$ be any limit process of subsequence of~$X^\varepsilon$.
Then, we have $X \in C([0,+\infty) ;\gamma)$ almost surely by (\ref
{th-curve1-1}).
Hence, taking $\varepsilon\downarrow0$ in~(\ref{SDE51}) with
replacing $X_x^\varepsilon$, $x$ and $\zeta_x^\varepsilon$ by
$X^\varepsilon$, $X^\varepsilon(0)$ and $\zeta^\varepsilon$, respectively,
\begin{eqnarray*}
\phi_1(X(t))&=&\phi_1(X(0)) + \int_0 ^t
\nabla\phi_1 (X(s)) \sigma(X(s))\,d\tilde{W}(s) \\
&&{}+ \int_0^t \nabla\phi_1 (X(s)) b(X(s))\,ds \\
&&{}+ \frac12 \sum_{i,j=1}^n \int_0 ^t a_{ij} (X(s))\,
\partial_i \partial_j \phi_1 (X(s))\,ds,
\end{eqnarray*}
where $\tilde W$ is an Wiener process.\vadjust{\goodbreak}

Noting that $\phi_1(X(\cdot))$ is a stochastic process on ${\mathbb
R}$ and $|\nabla\phi_1 (x) \sigma(x)|>0$ for $x\in \gamma$,
the law of $\phi_1(X(\cdot))$ is uniquely determined by this
equation; see Theorem~3.3 of Chapter IV in~\cite{IW}.
Applying It\^o's formula to $\gamma(\phi_1(X(t)))$ and noting that
$\gamma(\phi_1(X(\cdot)))=X(\cdot)$,
$\partial _i \phi _1 = \dot \gamma _i \circ \gamma ^{-1}$ on $\gamma$ for $i=1,2,\ldots ,N$
and $\partial _i \partial _j \phi _1 = (\dot \gamma _i \circ \gamma ^{-1}) (\ddot \gamma _j \circ \gamma ^{-1}) +
 (\dot \gamma _j \circ \gamma ^{-1}) (\ddot \gamma _i \circ \gamma ^{-1})$ on $\gamma$ for $i,j =1,2,\ldots ,N$,
we have that $X$ satisfies
(\ref{SDE10}); therefore, the first assertion holds.
The second assertion is obtained in a similar way.
\end{pf}

The argument above is also available in the case where the boundary
$\partial\Omega^\varepsilon$ carries a Neumann boundary condition,
for the generator of the process, in the following sense.
Consider a diffusion process $\widehat X^\varepsilon$ which is
associated with
\[
\frac12 \sum_{i,j=1}^n a_{ij}(x)\,\frac\partial{\partial x_i}\,\frac
\partial{\partial x_j} + \sum_{j=1}^n b_j(x)\,\frac\partial{\partial x_j}
\]
in $\Omega^\varepsilon$ and reflecting on $\partial\Omega^\varepsilon$.
Then, $\widehat X^\varepsilon$ can be expressed by the following equation:
%
\begin{equation}\label{SDE1-2}
\widehat X^\varepsilon(t)=\widehat X^\varepsilon(0) + \int_0 ^t
\sigma(\widehat X^\varepsilon(s))\,dW(s) + \int_0^t b(\widehat
X^\varepsilon(s))\,ds + \Phi^\varepsilon(\widehat X^\varepsilon)(t),\hspace*{-35pt}
\end{equation}
where $\Phi^\varepsilon$ is a singular drift which forces the
reflecting boundary condition on $\partial\Omega^\varepsilon$; see
\cite{SV2}.
Discussing this case in a similar way as above, we obtain the following theorem.

\begin{Th}\label{th-curve2}
Define a diffusion process $\widehat X$ by the solution of the
following equation:
%
\begin{eqnarray}\label
{SDE10-2}
 \widehat X(t)     &=&    \widehat X(0) + \int_0 ^t \dot
\gamma\circ\gamma^{-1} (\widehat X(s)) \langle\dot\gamma\circ
\gamma^{-1} (\widehat X(s)), \sigma(\widehat X(s))\,dW(s) \rangle\nonumber\\
&&{} + \int_0^t \dot\gamma\circ\gamma^{-1} (\widehat X(s)) \langle
\dot\gamma\circ\gamma^{-1} (\widehat X(s)), b(\widehat X(s))
\rangle \,ds\nonumber \\
&&{} + \frac 12 \int _0 ^t \ddot \gamma \circ \gamma ^{-1} (\widehat X(s))
 | \sigma (\widehat X(s))^T \dot \gamma \circ \gamma ^{-1} (\widehat X(s))| ^2 \,ds \\
&&{} + \int _0 ^t \dot \gamma \circ \gamma ^{-1} (\widehat X(s)) \langle \sigma (\widehat X(s))^T \ddot \gamma \circ \gamma ^{-1} (\widehat X(s)),\nonumber \\
&&\hspace*{100pt}{}  \sigma (\widehat X(s))^T \dot \gamma \circ \gamma ^{-1} (\widehat X(s))\rangle\, ds.\nonumber
\end{eqnarray}

If $\widehat X^\varepsilon(0)$ converges to a $\gamma$-valued random
variable $\widehat X(0)$ weakly, then the process $\widehat
X^\varepsilon$ converges weakly to $\widehat X$ in the sense of their
laws on $C([0,+\infty); {\mathbb R}^n)$ as $\varepsilon\downarrow0$.
Moreover, if pathwise uniqueness holds for~(\ref{SDE10-2}) and (\ref
{SDE1-2}) for all $\varepsilon>0$, and $\widehat X^\varepsilon(0)$
converges to a $\gamma$-valued random variable $\widehat X(0)$ almost
surely, then $\widehat X^\varepsilon$ converges to $\widehat X$ almost
surely, as $\varepsilon\downarrow0$.
\end{Th}

\begin{re}\label{re-curve}
In this section, the shape of tubes was taken to be cylindrical and the
``confining'' potential $U^\varepsilon$\vadjust{\goodbreak} has been defined by the scaling
of a fixed function~$U$.
However, neither the shape of the tubes nor the scaling property are essential.
If $U^\varepsilon$ is ``along $\gamma$'' (in the sense that the gradient
of~$U^\varepsilon$ is normal to the tangent of $\gamma$),
the same results hold.
In the case where~$U^\varepsilon$ is not along $\gamma$, some effect
of~$U^\varepsilon$ remains in the limit process; see~\cite{Sp,FS}.
\end{re}

\section{The case of $N$-spiders}\label{section spider}

In this section, we consider the shrinking of thin tubes to $N$-spider graphs.
The argument in this section is the main part of this article.
Consider an $n$-dimensional Euclidean space ${\mathbb R}^n$, let
$d(\cdot,\cdot)$ be the distance function in ${\mathbb R}^n$ and let
$O$ be the origin.
Let $\{e_i\} _{i=1}^N$ be $N$ different unit vectors in ${\mathbb R}^n$
and $ I_i :=\{ se_i \dvtx  s\in[0,\infty)\}$.
Consider an $N$-spider graph~$\Gamma$ defined by $ \Gamma:= \bigcup
_{i=1} ^N I_i $.
$\Gamma$ is also called an $N$-star graph.
Let $A$ be the set in ${\mathbb R}^n$ given by
\[
A:= \bigcup_{{i,j}\dvtx  {i \neq j}}  \{ x \in{\mathbb R}^n \dvtx  \langle
x,
e_i\rangle = \langle x, e_j\rangle  \} .
\]
For $x\in{\mathbb R}^n\setminus A$, let $\pi(x)$ be the nearest point
in $\Gamma$ from $x$.
Note that $\pi(x)$ is uniquely determined for all $x\in{\mathbb
R}^n\setminus A$.

Let $u_i$ be given similarly to $u$ in Section~\ref{section curve} for
$i=1,2,\ldots,N$ (so that~$u_i$ determines the potential acting in the
thin tube around $I_i$).
Let $c_i$ be a~positive number for $i=1,2,\ldots,N$,
\[
\kappa:= \max \bigl\{2 \sqrt2 c_i/\sqrt{1-\langle e_i,e_j\rangle} \dvtx
i,j=1,2,\ldots,N, i\neq j \bigr\}
\]
and $\kappa_0 \in(0,
\kappa)$.
$c_i$ has the interpretation of width of the tube around $I_i$.
Let~$U$ be a function on ${\mathbb R}^n$ with values in $[0,\infty]$,
and assume
\begin{eqnarray*}
  U(x) &= & u_i(c_i ^{-1}d(x,\Gamma)), \qquad  x \in\{
x\in{\mathbb R}^n\dvtx  \pi(x)\in I_i, d(x,I_i) <c_i, |x|\geq\kappa\} ,\\
U(x) &= & +\infty, \qquad  x \in\{ x\in{\mathbb R}^n\dvtx  \pi(x)\in I_i,
d(x,I_i) \geq c_i, |x| \geq\kappa\}, \\
U(x) &< & +\infty, \qquad  x\in\{ x\in{\mathbb R}^n\dvtx  |x|\leq\kappa_0\} ,
\end{eqnarray*}
$\Omega:= \{ x\dvtx  U(x)<\infty\}$ is a simply connected and unbounded
domain, $\partial\Omega$ is a~$C^2$-manifold and $U|_{\Omega}$ is a
$C^1$-function in $\Omega$.
This structure $\Omega$ is sometimes called a ``fattened'' $N$-spider.
In addition, we assume
\[
\lim _{m\rightarrow \infty} \langle -\nabla U(x_m), \nabla d(x_m , \partial \Omega)\rangle = +\infty \quad\mbox{and}\quad
-\lim_{m\rightarrow\infty}\frac{U(x_m)}{  \log  (d(x_m,\partial
\Omega))} = + \infty
\]
for any sequence $\{ x_m\}$ which converges to a point $x \in\partial
\Omega$.
Define domains $\{ \Omega_i \dvtx i=1,2,\ldots,N\}$ in ${\mathbb R}^n$ by
\[
\Omega_i := \{ x \in\Omega\setminus A\dvtx   \pi(x)\in I_i, |x|\geq
\kappa\}
\]
for $i=1,2,\ldots,N$.
Let $\Omega^\varepsilon:= \varepsilon\Omega$, $\Omega_i^\varepsilon
:= \varepsilon\Omega_i$, and $ U^\varepsilon(x)= U(\varepsilon
^{-1}x) $ for $x \in{\mathbb R}^n$ for all $\varepsilon>0$.
Note that $ U^\varepsilon(x) \in[0,+\infty)$ for $x\in\Omega
^\varepsilon$, $\partial U^\varepsilon$ is\vadjust{\goodbreak} a $C^2$-manifold, and
$U^\varepsilon|_{\Omega^\varepsilon}$ is a $C^1$-function on $\Omega
^\varepsilon$.
Consider a diffusion process $X^\varepsilon$ given by the following equation:
%
\begin{eqnarray}\label{SDE2}
X^\varepsilon(t)&=& X^\varepsilon(0) + \int_0 ^{t\wedge\zeta
^\varepsilon}\sigma(X^\varepsilon(s))\,dW(s) + \int_0^{t\wedge\zeta
^\varepsilon}b(X^\varepsilon(s))\,ds \nonumber
\\[-8pt]
\\[-8pt]
&& {}- \int_0^{t\wedge\zeta^\varepsilon}(\nabla U^\varepsilon
)(X^\varepsilon(s))\,ds,
\nonumber
\end{eqnarray}
where $X^\varepsilon(0)$ is an $\Omega^\varepsilon$-valued random
variable, $\zeta^\varepsilon$ is the first hitting time of
$X^\varepsilon$ at $\partial\Omega^\varepsilon$, $W$ is an
$n$-dimensional Wiener process, $ \sigma\in C_b({\mathbb R}^n;
{\mathbb R}^n\otimes{\mathbb R}^n) $ and $ b \in C_b({\mathbb R}^n;
{\mathbb R}^n)$.
Define a stochastic process $X^\varepsilon_x$ by the solution of (\ref
{SDE2}) with replacing $X^\varepsilon(0)$ by $x$, and $P_x^\varepsilon
$ by the law of $X^\varepsilon_x$ on $C([0,\infty); {\mathbb R}^n)$.
Let $ a(x) := \sigma(x) \sigma^T(x) $, and assume that $a$ is a
uniformly positive definite matrix.
Define a second-order elliptic differential operator $L$ on $\Omega
^\varepsilon$ by
\[
L:= \frac12 \sum_{i,j=1}^n a_{ij}(x)\,\frac{\partial}{\partial
x_i}\,\frac{\partial}{\partial x_j}+ \sum_{i=1}^n b_i(x)\,\frac
{\partial}{\partial x_i};
\]
then the generator of $X^\varepsilon$ is a closed extension of $
(L-\nabla U^\varepsilon\cdot\nabla)$ in $L^2(\Omega^\varepsilon
,dx)$ for any $\varepsilon>0$.
Since $a$ is a uniformly positive definite matrix, the process
$X^\varepsilon$ exists uniquely for all $\varepsilon>0$.

The following lemma implies that $X^\varepsilon$ does not exit from
$\Omega^\varepsilon$ almost surely.

\begin{lm}\label{lm2-1}
$\zeta^\varepsilon= +\infty$ almost surely for all $\varepsilon>0$.
\end{lm}

\begin{pf}
Locally, the discussion in the proof of Lemma~\ref{lm1-1} is available.
Hence, by using the strong Markov property of $X^\varepsilon$, we have
the assertion.
\end{pf}

Next we shall study the tightness of $\{ X^\varepsilon\dvtx  \varepsilon>0\}$.

\begin{lm}\label{lm-tight}
If the laws of $\{ X^\varepsilon(0) \dvtx  \varepsilon>0\}$ are tight,
then the laws of $\{ X^\varepsilon\dvtx   \varepsilon>0\}$ are also tight
in the sense of laws on $C([0,\infty) ;{\mathbb R}^n)$.
\end{lm}

\begin{pf}
In view of Theorem 2.1 in~\cite{FW} it is sufficient to show that for
any $\rho>0$ there exists a positive constant $C_ \rho$ such that for
all $y \in{\mathbb R}^n$ there exists a function $f_\rho^y$ on
${\mathbb R}^n$ which satisfies the following:
\begin{longlist}[(ii)]
\item[(i)] \hypertarget{tight1} $ f_\rho^y (y)=1 $, $f_\rho^y (x)=0$ for $|x-y|\geq\rho
$ and $0\leq f_\rho^y \leq1$.
\item[(ii)]\hypertarget{tight2}$ (f_\rho^y(X^\varepsilon(t))+ C_\rho t ; t\geq0)$
is a submartingale for sufficiently small $\varepsilon$.
\end{longlist}

Now we choose $f_\rho^y$ and $C_ \rho$ satisfying the conditions above.
Fix $\rho>0$, and take $\varepsilon_0 >0$ such that $\varepsilon_0<
\rho/{(16\kappa)}$.
When $y \in\overline{\Omega^{\varepsilon_0}}$ (where $\overline
{\Omega^{\varepsilon_0}}$ denotes the closure of $\Omega
^{\varepsilon_0}$ in ${\mathbb R}^n$) and $|y|>\rho/2$, choose
$f_\rho^y \in C^\infty({\mathbb R}^n)$ such that:
\begin{itemize}
\item$f_\rho^y(x)= f_\rho^y(\pi(x))$ for $x \in\Omega
^{\varepsilon_0}\setminus A$ and $f_\rho^y(x)=0$ for $|x-y|\geq\rho/4$;
\item$f_\rho^y(y)=1$, $0\leq f_\rho^y \leq1$, $\|\nabla f\|_\infty
\leq8/\rho$ and $\|\nabla^2f\|_\infty\leq64/\rho^2$.
\end{itemize}
Since $f_\rho^y(x)=0$ for $|x| \leq2\kappa\varepsilon_0$ and
$\nabla\pi(x) \nabla U^\varepsilon(x)=0$ for $|x| \geq2\kappa
\varepsilon_0$, it follows by It\^o's formula that
\[
f_\rho^y(X^\varepsilon(t)) - \int_0^t L f_\rho^y (X^\varepsilon(s))\,ds
\]
is a martingale for all $\varepsilon< \varepsilon_0$. Hence, choosing
$C_\rho$ larger than $ (8/\rho+ 64/{\rho^2})\times (\|\sigma\|_\infty
^2 /2+ \|b\|_\infty)$, conditions \hyperlink{tight1}{(i)} and \hyperlink{tight2}{(ii)} are satisfied for
$\varepsilon<\varepsilon_0$.\vspace*{1pt}

When $y \in\overline{\Omega^{\varepsilon_0}}$ and $|y| \leq\rho
/2$, choose $f_\rho^y \in C^\infty({\mathbb R}^n)$ such that:
\begin{itemize}
\item$f_\rho^y(x)= f_\rho^y(\pi(x))$ for $x \in\Omega
^{\varepsilon_0}\setminus A$, $f_\rho^y(x)=1$ for $|x| \leq\rho/4$,
and $f_\rho^y(x)=0$ for $|x-y|\geq\rho$;
\item $f_\rho ^y (y)=1$, $0\leq f_\rho^y \leq1$, $\|\nabla f\|_\infty\leq8/\rho$ and
$\|\nabla^2f\|_\infty\leq64/\rho^2$.
\end{itemize}
Here, note that $2\kappa\varepsilon\leq\rho/4$ for $\varepsilon
<\varepsilon_0$. Similarly to the case where $y \in\overline{\Omega
^{\varepsilon_0}}$ and $|y|>\rho/2$, one proves that conditions \hyperlink{tight1}{(i)}
and \hyperlink{tight2}{(ii)} are satisfied for $\varepsilon<\varepsilon_0$ with the same
$C_\rho$ as above.

When $y \notin\overline{\Omega^{\varepsilon_0}}$, choose $f_\rho
^y \in C^\infty({\mathbb R}^n)$ such that $f_\rho^y(y)=1$, $f_\rho
^y(x)=0$ for $x\in\overline{\Omega^{\varepsilon_0}}$, and $f_\rho
^y$ satisfies condition \hyperlink{tight1}{(i)} above.
Since $X^\varepsilon$ moves in $\Omega^\varepsilon$, $f_\rho
^y(X^\varepsilon(t))=0$ for all $t$ and $\varepsilon<\varepsilon_0$.

Thus, for all $\rho>0$, $\{ f_\rho^y\dvtx  y\in{\mathbb R}^n\}$ and $C_p$
are chosen in such a way that conditions \hyperlink{tight1}{(i)} and
\hyperlink{tight2}{(ii)} are satisfied.
\end{pf}

Now, we assume the tightness of $\{ X^\varepsilon(0)\dvtx  \varepsilon>0\}$.
By Lemma~\ref{lm-tight} we can choose a subsequence $\{ X^{\varepsilon
'}\dvtx  {\varepsilon'}>0\}$ of $\{ X^\varepsilon\dvtx  \varepsilon>0\}$ such
that the laws of its members converge weakly in the sense of laws on
$C([0,\infty); {\mathbb R}^n)$.
Define $X$~as the limit process of this subsequence, and to simplify~the
notation denote~the subsequence $\varepsilon'$ by $\varepsilon$ again.
From now on we fix $X$ as the limit process of $X^\varepsilon$.

For $w\in C([0,+\infty); {\mathbb R}^n)$, let $ \tilde T^c(w) := \inf
\{ t>0\dvtx  |w(t)| =c \}$ and $ T^c(w):= \inf\{ t>0 \dvtx  w(t)\notin A, |\pi
(w(t))|=c\} $ for $c>0$.

Theorem~\ref{th-curve} determines the behavior of $X$ on $\Gamma
\setminus O$.
Hence, to characterize~$X$, we need to determine the boundary condition
for $X$ at $O$.
Now we give some lemmas.
The following lemma implies that the edge which $X$ goes to, starting
from $O$, is independent of the edge which $X$ comes from.
Therefore, we obtain in particular that $X$ is a strong Markov process
on $\Gamma$.

\begin{lm}\label{lm-es2}
Let $\{ \delta(\varepsilon)\dvtx \varepsilon>0 \}$ be positive numbers
satisfying the condition that $\lim_{\varepsilon\downarrow
0}\varepsilon^{-1}\delta(\varepsilon)=+\infty$.
For $B \in{\mathscr B}({\mathbb R}^n)$ [${\mathscr B}({\mathbb R}^n)$
denoting the Borel subsets of ${\mathbb R}^n$],
\[
\sup \bigl\{  \bigl| P_x^\varepsilon \bigl( w\bigl(T ^{\delta(\varepsilon
)}\bigr) \in B  \bigr) - P_O^\varepsilon \bigl( w\bigl(T ^{\delta(\varepsilon
)}\bigr)\in B \bigr)  \bigr| \dvtx  x \in\Omega^\varepsilon,|x|\leq3\kappa
\varepsilon \bigr\}
\]
converges to $0$ as $\varepsilon\downarrow0$.
\end{lm}

\begin{pf}
Define a process ${\widehat X}_x^\varepsilon$ by the solution of the equation
%
\begin{eqnarray}\label{lm-es2-SDE}
  {\widehat X}_x^\varepsilon(t) &=& x + \int_0 ^t \sigma(\varepsilon
{\widehat X}_x^\varepsilon(s))\,d{\widehat W}(s) + \varepsilon\int_0^t
b(\varepsilon{\widehat X}_x^\varepsilon(s))\,ds\nonumber
\\[-8pt]
\\[-8pt]
&&{} - \int_0^t (\nabla
U)({\widehat X}_x^\varepsilon(s))\,ds
\nonumber
\end{eqnarray}
for $x\in\Omega$ and $\varepsilon>0$, where ${\widehat W}$ is an
$n$-dimensional Wiener process defined by $ {\widehat W}(t)=
\varepsilon^{-1}W(\varepsilon^2 t)$ for $t\in[0,\infty)$.
It is easy to see that the law of\vspace*{1pt} $( \widehat X ^\varepsilon_x (t)
\dvtx\break
t\geq0 )$ is equal to $(\varepsilon^{-1} X ^\varepsilon_{\varepsilon
x} (\varepsilon^2 t) \dvtx  t\geq0 )$ for $x \in\Omega$.
Letting ${\widehat P}_x^\varepsilon$ be the law of~${\widehat
X}_x^\varepsilon$ on $C([0,\infty); {\mathbb R}^n)$, we have
%
\begin{equation}\label{scaleP}
{\widehat P}_x^\varepsilon\bigl(w(t)\in dx\bigr) = P _{\varepsilon x}
^\varepsilon\bigl(\varepsilon^{-1} w(\varepsilon^2 t)\in dx\bigr)
\end{equation}
for $t\in[0,\infty)$, $x \in\Omega$ and $\varepsilon>0$.
By~(\ref{scaleP}), it is sufficient to show that
%
\begin{equation}\label{abc}
 \bigl| {\widehat P}_x^\varepsilon \bigl( w\bigl(T ^{\delta(\varepsilon
)/\varepsilon}\bigr) \in\varepsilon^{-1}B \bigr) - {\widehat
P}_O^\varepsilon\bigl ( w\bigl(T ^{\delta(\varepsilon)/\varepsilon}\bigr) \in
\varepsilon^{-1}B \bigr)  \bigr| \rightarrow0
\end{equation}
as $\varepsilon$ tends to $0$, uniformly in $x \in\{ y\in\Omega
\dvtx |y|\leq3 \kappa\}$.
Define stopping times
\begin{eqnarray*}
\tau_0 (w)&:=& \inf\{ t> 0 \dvtx  w(t)\notin A, |\pi(w(t))|> 3\kappa\}
,\\
\tilde\tau_{k} (w)&:=& \inf\{ t> \tau_{k-1} \dvtx  w(t)\notin A, |\pi
(w(t))|> 4 \kappa\},\qquad k\in{\mathbb N},\\
\tau_{k} (w)&:=& \inf\{ t> \tilde\tau_k \dvtx  w(t)\notin A, |\pi
(w(t))|< 3 \kappa\},\qquad k\in{\mathbb N},
\end{eqnarray*}
for $w \in C([0,\infty); {\mathbb R}^n)$. Note that $ |w(\tau_k)|
=3\kappa$ for $k=0,1,2,\ldots,$ and $ |w(\tilde\tau_k)| =4 \kappa
$ for $k=1,2,3,\ldots$ almost surely under $\widehat P^\varepsilon_x$
for $x\in\Omega$ and $|x|\leq2\kappa\varepsilon$.
Since $\triangle\pi(x)=0 $ and $\nabla\pi(x) \nabla U(x) =0$ for
$|x|\geq2\kappa$, It\^o's formula implies
%
\begin{eqnarray}\label{SDE-pi}
\pi({\widehat X}_x^\varepsilon(t)) &=& \pi(\tilde\tau_k ({\widehat
X}_x^\varepsilon)) + \int_{\tilde\tau_k ({\widehat X}_x^\varepsilon
)}^t \nabla\pi({\widehat X}_x^\varepsilon(s)) \sigma(\varepsilon
{\widehat X}^\varepsilon(s))\,d{\widehat W}(s) \nonumber
\\[-8pt]
\\[-8pt]
&&{} + \varepsilon\int_{\tilde\tau_k ({\widehat X}_x^\varepsilon
)}^t \nabla\pi({\widehat X}_x^\varepsilon(s)) b(\varepsilon
{\widehat X}^\varepsilon(s))\,ds
\nonumber
\end{eqnarray}
for $ t \in[\tilde\tau_k ({\widehat X}_x^\varepsilon),\tau_k
({\widehat X}_x^\varepsilon)]$, $ x\in\Omega$ and $|x|\leq
3\kappa\varepsilon$.
Since the diffusion coefficient of the one-dimensional process $|\pi
({\widehat X}_x^\varepsilon(t))|$ is uniformly elliptic, and
$T^{\delta(\varepsilon)/\varepsilon}$ diverges to infinity as
$\varepsilon\downarrow0$ almost surely under ${\widehat
P}_x^\varepsilon$, there exists a sequence $\{ \eta(\varepsilon)\}$
converging to $0$ as $\varepsilon\downarrow0$ such that
\[
\sup_{|x|=4\kappa}{\widehat P}_x^\varepsilon \bigl( T ^{\delta
(\varepsilon)/\varepsilon} < T^{3\kappa} \bigr) \leq\eta
(\varepsilon).
\]
On the other hand, since $\sigma\sigma^T$ is uniformly positive
definite, ${\widehat X}_x^\varepsilon$ hits $\{ x\in\Omega\dvtx\break  |x|<
\delta'\}$ with positive probability for all $x\in\Omega$,
$\varepsilon>0$, $\delta'>0$.
Hence, letting $\alpha(\varepsilon)$ be a sequence of positive
numbers such that $ \alpha(\varepsilon) \leq2\kappa$, and $\alpha
(\varepsilon)$ converges to $0$ as $\varepsilon\downarrow0$, we
obtain that
\[
p(\varepsilon):=\inf_{|x|=3\kappa}{\widehat P}_x^\varepsilon \bigl(
\tilde T ^{\alpha(\varepsilon)} < T^{4\kappa} \bigr) >0\vadjust{\goodbreak}
\]
for all $\varepsilon>0$, and that $p(\varepsilon)$ converges to $0$
as $\varepsilon\downarrow0$.
Moreover, we have
\begin{eqnarray*}
&&\hspace*{-5pt} {\widehat P}_x^\varepsilon \bigl( T^{\delta(\varepsilon
)/\varepsilon} < \tilde T^{\alpha(\varepsilon)} \bigr) \\
&&\hspace*{-5pt} \qquad  = \sum_{k=1}^{\infty} {\widehat P}_x^\varepsilon \bigl( T ^{\delta
(\varepsilon)/\varepsilon} < \tau_k, \tilde\tau_k<\tilde T^{\alpha
(\varepsilon)} \bigr) \\
&&\hspace*{-5pt} \qquad  = \sum_{k=1}^{\infty} \int_{\{ x_1\in\Omega\dvtx  |\pi(x_1)|=3\kappa\}
}\int_{\{ y_1\in\Omega\dvtx  |\pi(y_1)|=4\kappa\}} \\
&&\hspace*{-5pt} \qquad  \quad     \cdots\int_{\{ x_k\in\Omega\dvtx  |\pi(x_k)|=3\kappa\}} \int
_{\{ y_k\in\Omega\dvtx  |\pi(y_k)|=4\kappa\}} {\widehat P}_{y_k}^\varepsilon
 \bigl( T ^{\delta(\varepsilon)/\varepsilon} < T^{3\kappa} \bigr) \\
&&\hspace*{-5pt} \qquad  \hphantom{   \cdots\int_{\{ x_k\in\Omega\dvtx  |x_k|=3\kappa\}} \int
_{\{ y_k\in\Omega\dvtx  |y_k|=4\kappa\}}}
{} \times{\widehat P}_{x_k}^\varepsilon \bigl( w(T^{4\kappa})
\in dy_{k}, T^{4\kappa} < \tilde T^{\alpha(\varepsilon)} \bigr)\\
&&\hspace*{-5pt} \qquad  \hphantom{   \cdots\int_{\{ x_k\in\Omega\dvtx  |x_k|=3\kappa\}} \int
_{\{ y_k\in\Omega\dvtx  |y_k|=4\kappa\}}}
{} \times
{\widehat P}_{y_{k-1}}^\varepsilon \bigl( w(T^{3\kappa}) \in dx_{k}, T
^{\delta(\varepsilon)/\varepsilon} > T^{3\kappa} \bigr) \\
&&\hspace*{-5pt} \qquad  \hphantom{     \cdots\int_{\{ x_k\in\Omega\dvtx  |x_k|=3\kappa\}} \int
_{\{ y_k\in\Omega\dvtx  |y_k|=4\kappa\}}}
 {} \times\cdots\times{\widehat P}_{x_1}^\varepsilon \bigl(
w(T^{4\kappa}) \in dy_{1}, T^{4\kappa} < \tilde T^{\alpha
(\varepsilon)} \bigr)\\
&&\hspace*{-5pt} \qquad  \hphantom{   \cdots\int_{\{ x_k\in\Omega\dvtx  |x_k|=3\kappa\}} \int
_{\{ y_k\in\Omega\dvtx  |y_k|=4\kappa\}}}
 {} \times
 {\widehat P}_x^\varepsilon \bigl( w(T^{3\kappa
}) \in dx_1, T ^{\delta(\varepsilon)/\varepsilon} > T^{3\kappa
} \bigr) \\
&&\hspace*{-5pt} \qquad  \leq\eta(\varepsilon) \sum_{k=1}^{\infty}  \bigl( 1-p(\varepsilon
)  \bigr) ^k \\
&&\hspace*{-5pt} \qquad  = \frac{\eta(\varepsilon)(1-p(\varepsilon))}{p(\varepsilon)}.
\end{eqnarray*}
Hence, if $ \eta(\varepsilon) /p(\varepsilon)$ converges to $0$ as
$\varepsilon\downarrow0$, $ {\widehat P}_x^\varepsilon (
T^{\delta(\varepsilon)/\varepsilon} < \tilde T^{\alpha(\varepsilon
)} ) $ converges to $0$ as $\varepsilon\downarrow0$.
Now we choose $\alpha(\varepsilon)$ so that $ \eta(\varepsilon)
/p(\varepsilon) $ converges to $0$ as $\varepsilon\downarrow0$.
Then ${\widehat P}_x^\varepsilon ( T^{\delta (\varepsilon )/\varepsilon } < \tilde T^{\alpha (\varepsilon )})$ converges to $0$ as
$\varepsilon\downarrow0$.
Thus, for~(\ref{abc}), it is sufficient to prove that
%
\begin{equation}\label{abc2}
\sup_{|x|\leq\alpha(\varepsilon)} \bigl| {\widehat P}_x^\varepsilon
 \bigl( w\bigl(T ^{\delta(\varepsilon)/\varepsilon}\bigr) \in\varepsilon
^{-1}B \bigr) - {\widehat P}_O^\varepsilon \bigl( w\bigl(T ^{\delta
(\varepsilon)/\varepsilon}\bigr) \in\varepsilon^{-1}B \bigr)  \bigr|
\rightarrow0
\end{equation}
as $\varepsilon\downarrow0$.
To show this convergence, we use the coupling method.
Let $\sigma^l \in C_b^\infty({\mathbb R}^n; {\mathbb R}^n \otimes
{\mathbb R}^n)$ and $b^l \in C_b^\infty({\mathbb R}^n; {\mathbb R}^n
)$ for $l=1,2,\ldots,$ such that
\[
\lim _{l\rightarrow \infty } \sup _{|x|\leq M}|\sigma ^l(x) -\sigma (x)| =0 \hspace*{9pt}\mbox{and}\hspace*{9pt}
\lim _{l\rightarrow \infty } \sup _{|x|\leq M} |b^l(x) -b(x)| =0\hspace*{21pt}  \mbox{for }M>0.
\]
%
Let $x$ be fixed, and consider a pair of stochastic processes
$(\widetilde X_x^{\varepsilon,l}, \widetilde X_O^{\varepsilon,l})$
defined by
\begin{eqnarray*}
{\widetilde X}_x^{\varepsilon,l}(t) &=& x + \int_0 ^t \sigma
^l(\varepsilon{\widetilde X}_x^{\varepsilon,l}(s))\,d{\widehat W}(s) \\
&& {}+ \varepsilon\int_0^t b^l(\varepsilon{\widetilde
X}_x^{\varepsilon,l}(s))\,ds - \int_0^t (\nabla U)({\widetilde
X}_x^{\varepsilon,l}(s))\,ds,\\
{\widetilde X}_O^{\varepsilon,l}(t) &=& \int_0 ^t \sigma
^l(\varepsilon{\widetilde X}_O^{\varepsilon,l}(s))H^{\varepsilon
,l}({\widetilde X}_x^{\varepsilon,l}(s), {\widetilde X}_O^{\varepsilon
,l}(s))\,d{\widehat W}(s) \\
&&{} + \varepsilon\int_0^t b^l(\varepsilon{\widetilde
X}_O^{\varepsilon,l}(s))\,ds - \int_0^t (\nabla U)({\widetilde
X}_O^{\varepsilon,l}(s))\,ds,
\end{eqnarray*}
where
\[
H^{\varepsilon,l}(x_1,x_2):= I_n - \frac{2\sigma
^l(\varepsilon x_2)^{-1}(x_1-x_2)(x_1-x_2)^T(\sigma^l(\varepsilon x_2)^{-1})^T}{|\sigma
^l(\varepsilon x_2)^{-1}(x_1-x_2)|^2}
\]
for $x_1,x_2 \in{\mathbb R}^n$, and $I_n$ is the unit matrix.
Note that $(\widetilde X_x^{\varepsilon,l}, \widetilde
X_O^{\varepsilon,l})$ is uniquely determined because of the smoothness
of $\sigma^l$ and $b^l$.
We define
\begin{eqnarray*}
&  V(y):= |y|^{-1} y,\qquad\Phi^{\varepsilon,l}(x_1,x_2):= \sigma
^l(\varepsilon x_1)- \sigma^l(\varepsilon x_2)H^{\varepsilon,l}(x_1
,x_2) ,& \\
& \Psi^{\varepsilon,l}(x_1,x_2):= \varepsilon b^l(\varepsilon
x_1)-(\nabla U)(x_1)- \varepsilon b^l(\varepsilon x_2)+(\nabla U)(x_2)&
\end{eqnarray*}
for $y\in\{ z \in{\mathbb R}^n\dvtx  |z|\leq2\kappa_0\}$ and $x_1,x_2
\in\{ z \in{\mathbb R}^n\dvtx  |z|\leq\kappa_0\}$.
Similarly to the argument in Section 3 in~\cite{LR}, there exists a
positive constant $K$ such that $\inf_{l}\{ |\Phi^{\varepsilon
,l}(x_1,x_2)^T V(x_1-x_2)|^2 \} \geq K$ for $x_1,x_2
\in\{ z \in{\mathbb R}^n\dvtx  |z|\leq\kappa_0\}$ for small~$\varepsilon$.
By the equi-continuity of $\{\sigma^l\}$, we can choose $\rho\in
(0,2\kappa_0)$ satisfying
\begin{eqnarray*}
&& 2\langle x_1-x_2, \Psi^{\varepsilon,l}(x_1,x_2)\rangle+ \operatorname
{trace}(\Phi^{\varepsilon,l}(x_1,x_2)\Phi^{\varepsilon
,l}(x_1,x_2)^T) \\
 && \quad {}  - |\Phi^{\varepsilon,l}(x_1,x_2)^TV(x_1-x_2)|^2 \leq K/12
\end{eqnarray*}
for $|x_1-x_2| <\rho, |x_1| \leq\kappa_0,  |x_2| \leq\kappa_0,
l=1,2,3,\ldots$ (see~\cite{LR}).
For $\varepsilon'\in[0,\rho)$, define a stopping time ${\mathscr
T}^{\varepsilon'}$ on $C([0,\infty); {\mathbb R}^n) \times
C([0,\infty); {\mathbb R}^n)$ by
\begin{eqnarray*}
{\mathscr T}^{\varepsilon'}(w,w') &:=& \inf\{ t>0\dvtx  |w(t)-w'(t)|\notin(\varepsilon,\rho), \\
&&\hspace*{12pt}\hphantom{\inf\{}   |w(t)|\geq\kappa_0,    \mbox{ or }    |w'(t)|\geq\kappa_0
\}
\end{eqnarray*}
for $w,w'\in C([0,\infty); {\mathbb R}^n)$.
By It\^o's formula and the choice of $\rho$, we have
\begin{eqnarray*}
&&\hspace*{-5pt} \rho^{2/3} P \bigl(  | \widetilde X_x^{\varepsilon,l}({\mathscr
T}^{\varepsilon'}(\widetilde X_x^{\varepsilon,l},\widetilde
X_O^{\varepsilon,l}))-\widetilde X_O^{\varepsilon,l}({\mathscr
T}^{\varepsilon'}(\widetilde X_x^{\varepsilon,l},\widetilde
X_O^{\varepsilon,l})) | = \rho \bigr) \\
&&\hspace*{-5pt} \qquad  \leq E [  | \widetilde X_x^{\varepsilon,l}({\mathscr
T}^{\varepsilon'}(\widetilde X_x^{\varepsilon,l},\widetilde
X_O^{\varepsilon,l}))-\widetilde X_O^{\varepsilon,l}({\mathscr
T}^{\varepsilon'}(\widetilde X_x^{\varepsilon,l},\widetilde
X_O^{\varepsilon,l})) | ^{2/3} ] \\
&&\hspace*{-5pt} \qquad = |x|^{2/3} -\frac19 E \biggl[ \int_0^{{\mathscr T}^{\varepsilon
'}(\widetilde X_x^{\varepsilon,l}, \widetilde X_O^{\varepsilon,l})}
|\widetilde X_x^{\varepsilon,l}(s)-\widetilde X_O^{\varepsilon
,l}(s)|^{-4/3}   \\
&&\hspace*{-5pt} \qquad  \quad  \hphantom{|x|^{2/3} -\frac19 E \biggl[\int_0^{{\mathscr T}^{\varepsilon
'}(\widetilde X_x^{\varepsilon,l}, \widetilde X_O^{\varepsilon,l})}}   {} \times \bigl|
\Phi^{\varepsilon,l}(\widetilde X_x^{\varepsilon,l}(s), \widetilde
X_O^{\varepsilon,l}(s))^T\\
&&\hspace*{-5pt} \qquad  \quad  \hphantom{|x|^{2/3} -\frac19 E \biggl[\int_0^{{\mathscr T}^{\varepsilon
'}(\widetilde X_x^{\varepsilon,l}, \widetilde X_O^{\varepsilon,l})}{}\times\bigl|}   {} \times
 V\bigl(\widetilde
X_x^{\varepsilon,l}(s)-\widetilde X_O^{\varepsilon,l}(s)\bigr)  \bigr|
^2\,ds  \biggr] \\
&&\hspace*{-5pt} \qquad  \quad {} + \frac23 E\biggl [ \int_0^{{\mathscr T}^{\varepsilon
'}(\widetilde X_x^{\varepsilon,l}, \widetilde X_O^{\varepsilon,l})}
|\widetilde X_x^{\varepsilon,l}(s)-\widetilde X_O^{\varepsilon
,l}(s)|^{-4/3}   \\
&&\hspace*{-5pt} \qquad  \quad\hphantom{{}+\frac23 E\biggl [ \int_0^{{\mathscr T}^{\varepsilon
'}(\widetilde X_x^{\varepsilon,l}, \widetilde X_O^{\varepsilon,l})}}
{}\times \bigl\{
2\langle\widetilde X_x^{\varepsilon,l}(s)-\widetilde X_O^{\varepsilon
,l}(s), \Psi^{\varepsilon,l}(\widetilde X_x^{\varepsilon,l}(s),
\widetilde X_O^{\varepsilon,l}(s))\rangle  \\[-2pt]
&&\hspace*{-5pt}\hphantom{{}+\frac23 E\biggl [ \int_0^{{\mathscr T}^{\varepsilon
'}(\widetilde X_x^{\varepsilon,l}, \widetilde X_O^{\varepsilon,l})}{}\times\bigl\{}
 \qquad  \quad   {}+ \operatorname
{trace} \bigl( \Phi^{\varepsilon,l}(\widetilde X_x^{\varepsilon
,l}(s), \widetilde X_O^{\varepsilon,l}(s))\\[-2pt]
&&\hspace*{-5pt}\hphantom{{}+\frac23 E\biggl [ \int_0^{{\mathscr T}^{\varepsilon
'}(\widetilde X_x^{\varepsilon,l}, \widetilde X_O^{\varepsilon,l})}{}\times\bigl\{{}+\operatorname
{trace} \bigl(}
 \qquad  \quad
{}\times\Phi^{\varepsilon
,l}(\widetilde X_x^{\varepsilon,l}(s), \widetilde X_O^{\varepsilon
,l}(s))^T \bigr) \\[-2pt]
&&\hspace*{-5pt}\hphantom{{}+\frac23 E\biggl [ \int_0^{{\mathscr T}^{\varepsilon
'}(\widetilde X_x^{\varepsilon,l}, \widetilde X_O^{\varepsilon,l})}{}\times\bigl\{}
 \qquad  \quad   {} -  \bigl| \Phi
^{\varepsilon,l}   (\widetilde X_x^{\varepsilon,l}(s),
\widetilde X_O^{\varepsilon,l}(s))^T\\[-2pt]
&&\hspace*{9pt}\hspace*{-5pt}\hphantom{{}+\frac23 E\biggl [ \int_0^{{\mathscr T}^{\varepsilon
'}(\widetilde X_x^{\varepsilon,l}, \widetilde X_O^{\varepsilon,l})}{}\times\bigl\{}
 \qquad  \quad \hspace*{65pt}  {}\times
 V\bigl(\widetilde
X_x^{\varepsilon,l}(s)-\widetilde X_O^{\varepsilon,l}(s)\bigr)  \bigr|
^2 \bigr\}\,ds  \biggr] \\[-2pt]
&&\hspace*{-5pt} \qquad \leq|x|^{2/3} - \frac K{18} E \biggl[ \int_0^{{\mathscr
T}^{\varepsilon'}(\widetilde X_x^{\varepsilon,l}(s), \widetilde
X_O^{\varepsilon,l}(s))} |\widetilde X_x^{\varepsilon
,l}(s)-\widetilde X_O^{\varepsilon,l}(s)|^{-4/3}\,ds \biggr] \\[-2pt]
&&\hspace*{-5pt} \qquad \leq|x|^{2/3} - \frac{K}{18\rho^{4/3}} E [ {\mathscr
T}^{\varepsilon'}(\widetilde X_x^{\varepsilon,l}, \widetilde
X_O^{\varepsilon,l}) ].
\end{eqnarray*}
Hence, letting $\varepsilon' \downarrow0$, we have the following two
estimates:
%
\begin{eqnarray} \label{corr1}
&& P \bigl(  | \widetilde X_x^{\varepsilon,l}({\mathscr
T}^{0}(\widetilde X_x^{\varepsilon,l},\widetilde X_O^{\varepsilon
,l}))-\widetilde X_O^{\varepsilon,l}({\mathscr T}^{0}(\widetilde
X_x^{\varepsilon,l},\widetilde X_O^{\varepsilon,l})) | = \rho
 \bigr), \nonumber
 \\[-9pt]
 \\[-9pt]
&& \qquad   \leq\rho^{-2/3} |x|^{2/3},  \nonumber\\[-27pt]
\nonumber
\end{eqnarray}
\begin{equation}
\label{corr2}
   E [ {\mathscr T}^{0}(\widetilde X_x^{\varepsilon
,l}, \widetilde X_O^{\varepsilon,l}) ] \leq \frac{18\rho
^{4/3}}{K}|x|^{2/3}.
\end{equation}
On the other hand, by It\^o's formula,
\begin{eqnarray*}
&& E \biggl[  \biggl\{  \biggl(  | \widetilde X_x^{\varepsilon
,l}({\mathscr T}^{0}(\widetilde X_x^{\varepsilon,l},\widetilde
X_O^{\varepsilon,l})) -x  | - \frac{\kappa_0}{2}  \biggr) _+
 \biggr\} ^2 \biggr]\\[-2pt]
&& \qquad = E \biggl[ \int_0^{{\mathscr T}^{0}(\widetilde X_x^{\varepsilon
,l},\widetilde X_O^{\varepsilon,l})}  \biggl(  | \widetilde
X_x^{\varepsilon,l}(s) -x  | - \frac{\kappa_0}{2}  \biggr) _+
|\widetilde X_x^{\varepsilon ,l}(s) -x|^{-1}
{\mathbb I}_{\{ |\widetilde X_x^{\varepsilon,l}(s) -x|\geq\kappa_0/2
\}}  \\[-2pt]
&& \qquad  \quad
\hphantom{E \biggl[ \int_0^{{\mathscr T}^{0}(\widetilde X_x^{\varepsilon
,l},\widetilde X_O^{\varepsilon,l})} }
 {}\times \biggl\{   2\langle \widetilde X_x^{\varepsilon ,l}(s) -x, \varepsilon b^l(\varepsilon \widetilde X_x^{\varepsilon ,l}(s))
-\nabla U(\widetilde X_x^{\varepsilon ,l}(s)) \rangle \\
&&\hphantom{E \biggl[ \int_0^{{\mathscr T}^{0}(\widetilde X_x^{\varepsilon[-2pt]
,l},\widetilde X_O^{\varepsilon,l})}{}\times\biggl\{ }
 \qquad  \quad {}   + \operatorname{trace}
[ \sigma^l(\varepsilon\widetilde X_x^{\varepsilon,l}(s))\sigma
^l(\varepsilon\widetilde X_x^{\varepsilon,l}(s))^T ]   \\[-2pt]
&&\hspace*{54.5pt}\hphantom{E \biggl[ \int_0^{{\mathscr T}^{0}(\widetilde X_x^{\varepsilon
,l},\widetilde X_O^{\varepsilon,l})}{}\times\biggl\{ }
 \qquad  \quad {}   - \biggl| \sigma
^l(\varepsilon\widetilde X_x^{\varepsilon,l}(s))^T \frac{\widetilde
X_x^{\varepsilon,l}(s) -x}{ | \widetilde X_x^{\varepsilon,l}(s)
-x  |} \biggr| ^2  \biggr\} \\[-2pt]
&&\hphantom{E \biggl[ \int_0^{{\mathscr T}^{0}(\widetilde X_x^{\varepsilon
,l},\widetilde X_O^{\varepsilon,l})} }
 \qquad  \quad       {}+ {\mathbb
I}_{\{ |\widetilde X_x^{\varepsilon,l}(s) -x|\geq\kappa_0/2 \}}
\biggl| \sigma^l(\varepsilon\widetilde X_x^{\varepsilon,l}(s))^T \frac
{\widetilde X_x^{\varepsilon,l}(s) -x}{ | \widetilde
X_x^{\varepsilon,l}(s) -x  |} \biggr| ^2\,ds \biggr] \\[-2pt]
&& \qquad  \leq C E [ {\mathscr T}^{0}(\widetilde X_x^{\varepsilon,l},
\widetilde X_O^{\varepsilon,l}) ],
\end{eqnarray*}
where $z_+ := \max\{ 0,z\}$ for $z\in\mathbb R$ and $C$ is a positive
constant independent of $l$ and~$x$.
This inequality together with~(\ref{corr2}) implies
%
\begin{equation}\label{corr3}
P \bigl(  | \widetilde X_x^{\varepsilon,l}({\mathscr
T}^{0}(\widetilde X_x^{\varepsilon,l},\widetilde X_O^{\varepsilon
,l}))-x | = \kappa_0  \bigr) \leq\frac{72C\rho^{4/3}}{\kappa
_0^2K} |x|^{2/3}.\vadjust{\goodbreak}
\end{equation}
Similarly, we have
%
\begin{equation}
P \bigl(  | \widetilde X_O^{\varepsilon,l}({\mathscr
T}^{0}(\widetilde X_x^{\varepsilon,l},\widetilde X_O^{\varepsilon
,l})) | = \kappa_0  \bigr) \leq\frac{72C'\rho^{4/3}}{\kappa
_0^2K} |x|^{2/3}, \label{corr4}
\end{equation}
where $C'$ is a positive constant.
Noting that $\widetilde X_x^{\varepsilon,l}$ and $\widetilde
X_O^{\varepsilon,l}$ converge to $\widehat X_x^{\varepsilon}$ and
$\widehat X_O^{\varepsilon}$ in law as $l\rightarrow+\infty$,
respectively, for each $\varepsilon$, by the coupling inequality (see
\cite{LR}) we have
\begin{eqnarray*}
&& \sup_{|x|\leq\alpha(\varepsilon)} \bigl| {\widehat
P}_x^\varepsilon \bigl( w\bigl(T ^{\delta(\varepsilon)/\varepsilon}\bigr) \in
\varepsilon^{-1}B \bigr) - {\widehat P}_O^\varepsilon \bigl( w\bigl(T
^{\delta(\varepsilon)/\varepsilon}\bigr) \in\varepsilon^{-1}B \bigr)
 \bigr| \\[-1pt]
&& \qquad  \leq\sup_{|x|\leq\alpha(\varepsilon)} \sup_{l}  \bigl| P \bigl(
\widetilde X_x^{\varepsilon,l}\bigl(T ^{\delta(\varepsilon)/\varepsilon
}\bigr) \in\varepsilon^{-1}B \bigr) - P \bigl( \widetilde X_O^{\varepsilon
,l}\bigl(T ^{\delta(\varepsilon)/\varepsilon}\bigr) \in\varepsilon
^{-1}B \bigr)  \bigr| \\[-1pt]
&& \qquad  \leq\sup_{|x|\leq\alpha(\varepsilon)} \sup_{l} P \bigl(  |
\widetilde X_x^{\varepsilon,l}({\mathscr T}^{0}(\widetilde
X_x^{\varepsilon,l},\widetilde X_O^{\varepsilon,l}))-\widetilde
X_O^{\varepsilon,l}({\mathscr T}^{0}(\widetilde X_x^{\varepsilon
,l},\widetilde X_O^{\varepsilon,l})) |\neq 0  \bigr) \\[-1pt]
&& \qquad  \leq\sup_{|x|\leq\alpha(\varepsilon)} \sup_{l} P \bigl(  |
\widetilde X_x^{\varepsilon,l}({\mathscr T}^{0}(\widetilde
X_x^{\varepsilon,l},\widetilde X_O^{\varepsilon,l}))-\widetilde
X_O^{\varepsilon,l}({\mathscr T}^{0}(\widetilde X_x^{\varepsilon
,l},\widetilde X_O^{\varepsilon,l})) | =\rho \bigr) \\[-1pt]
&& \qquad  \quad {}  + \sup_{|x|\leq\alpha(\varepsilon)} \sup_{l} P \bigl(
| \widetilde X_x^{\varepsilon,l}({\mathscr T}^{0}(\widetilde
X_x^{\varepsilon,l},\widetilde X_O^{\varepsilon,l}))-x | = \kappa
_0  \bigr)\\[-1pt]
&& \qquad  \quad {}  + \sup_{|x|\leq\alpha(\varepsilon)} \sup_{l} P\bigl (
| \widetilde X_O^{\varepsilon,l}({\mathscr T}^{0}(\widetilde
X_x^{\varepsilon,l},\widetilde X_O^{\varepsilon,l})) | = \kappa
_0  \bigr) .
\end{eqnarray*}
This inequality, together with~(\ref{corr1}),~(\ref{corr3}) and (\ref
{corr4}) yields~(\ref{abc2}).
\end{pf}

The next lemma implies that $O$ is not absorbing for $X$.

\begin{lm}\label{lm-es1}
\[
\int_0^t E \bigl[ {\mathbb I}_{\{ x\dvtx  |x|\leq\delta'\}}(X(s)) \bigr]\,ds =O(\delta')
\]
as $\delta'\downarrow0$, for all $t\geq0$.
\end{lm}

\begin{pf}
To simplify the notation, let $X^\varepsilon(0) = x^\varepsilon\in
\Omega^\varepsilon$.
It is sufficient to show that
\[
\int_0^t E \bigl[ {\mathbb I}_{\{ x\dvtx  |\pi(x)|\leq\delta'\}
}(X(s)) \bigr]\,ds =O(\delta')
\]
as $\delta'\downarrow0$.
By Fatou's lemma, we have
%
\begin{eqnarray} \label{lm-es1-2}
&& \int_0^t E \bigl[ {\mathbb I}_{\{ x\dvtx  |\pi(x)|\leq\delta'\}
}(X(s)) \bigr]\,ds\nonumber\\
&& \qquad  \leq\liminf_{\varepsilon\downarrow0} \int_0^t E\bigl [ {\mathbb
I}_{\{ x\dvtx  |\pi(x)|\leq3\kappa\varepsilon\}}(X^\varepsilon(s))
\bigr]\,ds \\
&& \qquad  \quad {}  + \liminf_{\varepsilon\downarrow0} \int_0^t E \bigl[
{\mathbb I}_{\{ x\dvtx  3\kappa\varepsilon\leq|\pi(x)|\leq\delta'\}
}(X^\varepsilon(s)) \bigr]\,ds. \nonumber
\end{eqnarray}
To show that the second term is $O(\delta')$ as $\delta'\downarrow
0$, let $f$ be a continuous function on $\mathbb R$ such that $
{\mathbb I}_{\{ x \in{\mathbb R}\dvtx  3\kappa\varepsilon \leq x \leq
\delta'\}} \leq f \leq{\mathbb I}_{\{ x \in{\mathbb R}\dvtx  2\kappa
\varepsilon\leq x \leq2\delta'\}}$ and $ F(x) := \int_0^x \int
_0^y f(z)\,dz \,dy$.
Noting that $ \pi(x) = \langle e_i,x\rangle e_i$ for $x \in\Omega
_i$ and $i=1,2,\ldots,N$, we have $\nabla\pi(x) \pi(x) = \pi(x)$
for $x \in\Omega$ such that $|x|\geq2\kappa$.
Since $\nabla\pi(x) \nabla U ^\varepsilon(x)=0$ and $\triangle\pi
(x)=0$ for $x\in\Omega^\varepsilon$ such that $|x| \geq2\kappa
\varepsilon$, we have
\begin{eqnarray*}
&& E[ F(|\pi(X^\varepsilon(t))|) ] - F(|\pi(x^\varepsilon)|)\\
&& \qquad = \frac12 \int_0^t E \biggl[ f(|\pi(X^\varepsilon(s))|)  \biggl|
\sigma(X^\varepsilon(s)) ^T \frac{\pi(X^\varepsilon(s))}{|\pi
(X^\varepsilon(s))|}  \biggr| ^2  \biggr]\,ds\\
&& \qquad  \quad {} + \int_0^t E\biggl [ F'(|\pi(X^\varepsilon(s))|)  \biggl\langle
\frac{\pi(X^\varepsilon(s))}{|\pi(X^\varepsilon(s))|},
b(X^\varepsilon(s))  \biggr\rangle \biggr]\,ds.
\end{eqnarray*}
It is easy to see that $ E [ |X^\varepsilon(t)|^2  ] $ is
dominated uniformly in $\varepsilon>0$.
Moreover, it holds that $ 0 \leq F' \leq2\delta'$ and $ 0 \leq
F(x) \leq2\delta'x$ for $x \in{\mathbb R}_+$.
Thus, by uniform ellipticity of $a=\sigma\sigma^T$, we have the
following estimate:
\[
\int_0^t E \bigl[{\mathbb I}_{\{ x \in{\mathbb R}\dvtx  3\kappa
\varepsilon \leq x \leq\delta'\}}(|\pi(X^\varepsilon(s))|) \bigr]\,ds \leq C\delta'
\]
for some constant $C$.
Hence,
\[
\liminf_{\varepsilon\downarrow0} \int_0^t E\bigl [ {\mathbb I}_{\{
x\dvtx  3\kappa\varepsilon\leq|\pi(x)|\leq\delta'\}}(X^\varepsilon
(s)) \bigr]\,ds =O(\delta')
\]
as $\delta'\downarrow0$.
This yields that the second term of~(\ref{lm-es1-2}) is equal to
$O(\delta')$ as $\delta' \downarrow0$.

The proof is finished by showing that
%
\begin{equation}\label{lm-es1-1}
\int_0^t E \bigl[ {\mathbb I}_{\{ x\dvtx  |\pi(x)|\leq3\kappa\varepsilon
\}}(X^\varepsilon(s)) \bigr]\,ds = O(\varepsilon)
\end{equation}
as $\varepsilon\downarrow0$.
Define stopping times $ \{ \tau_{k} ^\varepsilon, \tilde\tau_{k}
^\varepsilon\} $ by
\begin{eqnarray*}
\tau_0 ^\varepsilon(w)&:=& 0,\\
\tilde\tau_{k} ^\varepsilon(w)&:=& \inf\{ u> \tau_{k-1}^\varepsilon
(w) \dvtx  |\pi(w(u))|> 4\kappa\varepsilon\},\qquad k\in{\mathbb N},\\
\tau_{k} ^\varepsilon(w)&:=& \inf\{ u> \tilde\tau_k ^\varepsilon
(w) \dvtx  |\pi(w(u))|< 3\kappa\varepsilon\},\qquad k\in{\mathbb N},
\end{eqnarray*}
for $w\in C([0,\infty); {\mathbb R}^n)$.
Then,
\begin{eqnarray*}
&&\int_0^t E \bigl[ {\mathbb I}_{\{ x\dvtx  |\pi(x)|\leq3\kappa
\varepsilon\}}(X^\varepsilon(s)) \bigr]\,ds \\
&& \qquad  \leq\sum_{k=1}^\infty\int \biggl( \int T^{4\kappa\varepsilon}
(w)P_x^\varepsilon(dw)  \biggr) P_{x^\varepsilon} ^\varepsilon\bigl(w(\tau
_{k}^\varepsilon)\in dx, \tau_{k} ^\varepsilon \leq t\bigr) \\
&& \qquad  \leq\sup_{x\in\{ y\in\Omega\dvtx  |\pi(y)|=3\kappa\varepsilon\}
}\biggl ( \int T^{4\kappa\varepsilon} (w) P_x^\varepsilon(dw)  \biggr)
\sum_{k=1}^\infty P_{x^\varepsilon} ^\varepsilon(\tau_{k}
^\varepsilon \leq t) .
\end{eqnarray*}
By using the notation in the proof of Lemma~\ref{lm-es2}, we have
\[
\sup_{x\in\{ y\in\Omega\dvtx  |\pi(y)|=3\kappa\varepsilon\}} \int
T^{4\kappa\varepsilon} (w) P_x^\varepsilon(dw) = \varepsilon^2 \sup
_{x\in\{ y\in\Omega\dvtx  |\pi(y)|=3\kappa\}} \int T^{4\kappa} (w)
{\widehat P}_x^\varepsilon(dw) .
\]
It is easy to see that
\[
\sup_{\varepsilon>0} \sup_{x\in\{ y\in\Omega\dvtx  |\pi(y)|=3\kappa\}
} \int T^{4\kappa} (w) {\widehat P}_x^\varepsilon(dw)< +\infty.
\]
Hence, for~(\ref{lm-es1-1}), it is sufficient to show that
%
\begin{equation}\label{lm-es1-3}
\sum_{k=1}^\infty P_{x^\varepsilon}^\varepsilon(\tilde\tau_{k}
^\varepsilon \leq t) \leq C\varepsilon^{-1}
\end{equation}
for some constant $C$.
For $w\in C([0,\infty); \Omega)$, let ${\mathcal N}_t(w)$ be the
number of transitions of $w$ from the set $ \{ x\in\Omega
^\varepsilon\dvtx  |\pi(x)|=3\kappa\} $ to the set $ \{ x\in\Omega
^\varepsilon\dvtx  |\pi(x)|= 4\kappa\} $ during the time interval $[0,t]$.
Then,
%
\begin{equation}\label{lm-es1-4}
\sum_{k=1}^\infty P_{x^\varepsilon}^\varepsilon(\tilde\tau_{k}
^\varepsilon \leq t)= \int{\mathcal N}_{\varepsilon^{-2}t}(w)
{\widehat P}_{\varepsilon^{-1}x^\varepsilon}^\varepsilon(dw).
\end{equation}

Take $f\in C^\infty ([0,\infty ))$ such that
$f\geq 0$, $0\leq f'\leq 1$, $f''\geq 0$, $\operatorname{supp} f'' \subset [2\kappa ,3\kappa ]$, $f(x)=0$ for $x\leq 2\kappa$
and $f(3\kappa ) < f(4\kappa )$.
Define $\widehat Y^{\varepsilon ,i}_x$ by
\[
\widehat Y^{\varepsilon ,i}_x (t) := f(\langle e_i, \widehat X_x^\varepsilon (t)\rangle {\mathbb I}_{\Omega _i}( \widehat X_x^\varepsilon (t)))
\]
for $x\in \Omega$ and $i=1,2,\ldots ,N$.
Since $\langle e_i, \nabla U (x) \rangle =0$ for $x \in \{ \Omega _i \dvtx |x| \geq 2\kappa \varepsilon\}$, by It\^o's formula we have
\begin{eqnarray*}
\widehat Y_x^{\varepsilon ,i}(t) &=& f(\langle e_i, x\rangle {\mathbb I}_{\Omega _i}(x)) \\
&&{} + \int _0 ^t f'(\langle e_i, \widehat X_x^\varepsilon (s)\rangle {\mathbb I}_{\Omega _i}
( \widehat X_x^\varepsilon (s))) \langle e_i, \sigma (\varepsilon \widehat X_x^\varepsilon (s))\,d{\widehat W}(s)\rangle \\
&&{} + \varepsilon \int _0^t f'(\langle e_i, \widehat X_x^\varepsilon (s)\rangle {\mathbb I}_{\Omega _i}
( \widehat X_x^\varepsilon (s))) \langle e_i, b(\varepsilon \widehat X_x^\varepsilon (s))\rangle\, ds \\
&&{} + \frac 12 \int _0^t f''(\langle e_i, \widehat X_x^\varepsilon (s)\rangle {\mathbb I}_{\Omega _i}
( \widehat X_x^\varepsilon (s))) |\sigma (\varepsilon \widehat X_x^\varepsilon (s))^Te_i|^2\,ds.
\end{eqnarray*}
It is clear that
\[
E[ {\mathcal N}_{\varepsilon ^{-2}t}(\widehat X^\varepsilon _{\varepsilon ^{-1}x^\varepsilon})]
\leq \sum _{i=1}^N \sup _{x: |\pi (x)|\leq 4\kappa }E[ \widetilde{\mathcal N}_{\varepsilon ^{-2}t}
(\widehat Y^{\varepsilon ,i} _x)] ,
\]
where $\widetilde{\mathcal N}_t(w)$ is the number of up-crossing of $w$ for the interval
$[f(3\kappa ),f(4\kappa )]$ during the time interval $[0,t]$.
Hence, by~(\ref{lm-es1-3}) and~(\ref{lm-es1-4}), it is sufficient to show that
\begin{equation}\label{eq3.15}
\sup _{x: |\pi (x)|\leq 4\kappa }E[ \widetilde{\mathcal N}_{\varepsilon ^{-2}t}
(\widehat Y^{\varepsilon ,i} _x)] \leq C\varepsilon ^{-1}
\end{equation}
with a constant $C$ for all $i=1,2,\ldots ,N$.
Let $i$ be fixed and $m\in {\mathbb N}$.
Define~$\tau _k$ and $\tilde \tau _k$ by
\begin{eqnarray*}
\tilde \tau _0 &:=& 0, \\
\tau _0 &:=& \inf\{ u>0 \dvtx {\widehat Y_x^{\varepsilon ,i}}(u) \leq f(3\kappa ) \} ,\\
\tilde \tau _{k} &:=& \inf\{ u> \tau _{k-1} \dvtx {\widehat Y_x^{\varepsilon ,i}}(u) \geq f(4\kappa )\}, \qquad k\in {\mathbb N},\\
\tau _{k} &:=& \inf\{ u> \tilde \tau _k \dvtx {\widehat Y_x^{\varepsilon ,i}}(u) \leq f(3\kappa ) \},\qquad  k\in {\mathbb N}.
\end{eqnarray*}
Then,
\begin{eqnarray*}
&&   E[ \widehat Y_x^{\varepsilon ,i}(t\wedge \tilde \tau _m )]
-E[ \widehat Y_x^{\varepsilon ,i}(t\wedge \tau _0 )] \\
&& \qquad  = \sum _{k=1}^m E[ \widehat Y_x^{\varepsilon ,i}(\tilde \tau _k\wedge t)
 - \widehat Y_x^{\varepsilon ,i}(\tau _{k-1}\wedge t) ]\\
 &&\qquad\quad{}
  + \sum _{k=1}^{m-1} E[ \widehat Y_x^{\varepsilon ,i}(\tau _k\wedge t) - \widehat Y_x^{\varepsilon ,i}(\tilde \tau _k\wedge t) ] \\
&&\qquad   = \sum _{k=1}^m E[ \widehat Y_x^{\varepsilon ,i}(\tilde \tau _k\wedge t)
 -\widehat Y_x^{\varepsilon ,i}(\tau _{k-1}\wedge t) ] \\
&& \qquad\quad{}  + \varepsilon \sum _{k=1}^{m-1}
 E\biggl[ \int _{\tilde \tau _k\wedge t}^{\tau _k\wedge t} f'(\langle e_i, \widehat X_x^\varepsilon (s)\rangle
  {\mathbb I}_{\Omega _i}( \widehat X_x^\varepsilon (s))) \langle e_i, b(\varepsilon \widehat X_x^\varepsilon (s))\rangle\,  ds \biggr] \\
&& \qquad\quad{}    + \frac 12 \sum _{k=1}^{m-1}
 E\biggl[\int _{\tilde \tau _k\wedge t}^{\tau _k\wedge t} f''(\langle e_i,
 \widehat X_x^\varepsilon (s)\rangle {\mathbb I}_{\Omega _i}( \widehat X_x^\varepsilon (s)))
 |\sigma (\varepsilon \widehat X_x^\varepsilon (s))^Te_i|^2\,ds\biggr] .
\end{eqnarray*}
Since $f''\geq 0$, we have
\begin{eqnarray}\label{eq3.16}
&&\Biggl| \sum _{k=1}^m E[ \widehat Y_x^{\varepsilon ,i}
(\tilde \tau _k\wedge t) - \widehat Y_x^{\varepsilon ,i}(\tau _{k-1}\wedge t) ]
\Biggr|
\nonumber
\\[-8pt]
\\[-8pt]
\nonumber
 &&\qquad\leq E[ \widehat Y_x^{\varepsilon ,i}(t\wedge \tilde \tau _m )] + E[ \widehat Y_x^{\varepsilon ,i}(t\wedge \tau _0 )] + C_1 \varepsilon t
\end{eqnarray}
with a positive constant $C_1$.
Let
\begin{eqnarray*}
\tilde \tau _* &:=& \max \{ \tilde \tau _k \dvtx \tilde \tau _k \leq t, k=1,2,3,\ldots \} ,\\
\tau _* &:=& \max \{ \tau _k \dvtx \tau _k \leq t, k=1,2,3,\ldots \} , \\
M(t) &:=& \int _0 ^t f'(\langle e_i, \widehat X_x^\varepsilon (s)\rangle {\mathbb I}_{\Omega _i}
( \widehat X_x^\varepsilon (s))) \langle e_i, \sigma (\varepsilon \widehat X_x^\varepsilon (s))\,d{\widehat W}(s)\rangle .
\end{eqnarray*}
When $\tilde \tau _* \leq \tau _*$, $ \widehat Y_x^{\varepsilon ,i}(t) \leq f(4\kappa )$.
When $\tau _* \leq \tilde \tau _*$, $f''(\langle e_i, \widehat X_x^\varepsilon (s)\rangle {\mathbb I}_{\Omega _i}
( \widehat X_x^\varepsilon (s))) = 0$ for $s\in [\tilde \tau _* ,t]$. Thus, we have
\begin{eqnarray*}
\widehat Y_x^{\varepsilon ,i}(t) &=& \widehat Y_x^{\varepsilon ,i}(\tilde \tau _* ) + M(t) - M(\tilde \tau _* ) \\
&&{}+ \varepsilon \int _{\tilde \tau _*}^t f'(\langle e_i, \widehat X_x^\varepsilon (s)\rangle {\mathbb I}_{\Omega _i}
( \widehat X_x^\varepsilon (s))) \langle e_i, b(\varepsilon \widehat X_x^\varepsilon (s))\rangle \,ds.
\end{eqnarray*}
Hence,
\[
\widehat Y_x^{\varepsilon ,i}(t) \leq f(4\kappa ) + 2 \sup _{0\leq s \leq t} |M(s)| + C_2 \varepsilon t
\]
for $|x|\leq 4\kappa$ with a constant $C_2$.
By the Burkholder--Davis--Gundy inequality we have
\[
E\Bigl[ \sup _{0\leq s \leq t} \widehat Y_x^{\varepsilon ,i}(s)\Bigr] \leq f(4\kappa ) + 2C_3\sqrt t + C_2 \varepsilon t
\]
for $|x|\leq 4\kappa$ with a constant $C_3$.
Thus, letting $m\rightarrow +\infty$ on (3.16), we have for $|x|\leq 4\kappa$
\[
\bigl(f(4\kappa )-f(3\kappa ) \bigr) E\Bigl[ \widetilde{\mathcal N}_t(\widehat Y^{\varepsilon ,i} _x)\Bigr]
\leq 2 f(4\kappa ) +4C_3\sqrt t + (C_1 + 2C_2) \varepsilon t .
\]
Therefore, replacing $t$ by $\varepsilon ^{-2}t$, (3.15) is obtained.
\end{pf}


The lemmas above yield that the boundary condition at $O$ is a weighted
Kirchhoff boundary condition.
Hence, the next step is to determine the weights associated with the edges.
Let $Y_x^\varepsilon$ be a diffusion process defined by the solution of
the following stochastic differential equation:
%
\begin{equation}\label{SDEY}
Y_x^\varepsilon(t)=x+ \sigma(O) W(t) - \int_0^t (\nabla
U^\varepsilon)(Y_x^\varepsilon(s))\,ds.
\end{equation}
Note that $Y_x^\varepsilon$ is a special case of $X^\varepsilon$ with
the condition $X^\varepsilon(0)=x$, and $Y_x^\varepsilon$ does not hit
$\Omega^\varepsilon$ almost surely.
Denote the law of $Y_x^\varepsilon$ on $C([0,\infty);{\mathbb R}^n)$
by $Q_x^\varepsilon$.
It is easy to see that the law of $Y_x^\varepsilon$ is the same as that
of $\varepsilon Y_{\varepsilon^{-1}x}^1(\varepsilon^{-2})$.
By~(\ref{lm-es2-SDE}) one has that the law of $\widehat X_x^\varepsilon
$ converges to that of $Y_x^1$ as $\varepsilon\downarrow0$, and
therefore, the law of $X_x^\varepsilon$ and that of $Y_x^\varepsilon$
are getting closer as $\varepsilon\downarrow0$.
In particular, we have
\[
\lim_{\varepsilon\downarrow0}  \bigl| P_O^\varepsilon
\bigl(w(T^{c\varepsilon})\in\Omega_i^\varepsilon\bigr) - Q_O^\varepsilon
\bigl(w(T^{c\varepsilon})\in\Omega_i^\varepsilon\bigr) \bigr| =0
\]
for all $c>0$ and $i=1,2,\ldots,N$.
Since this holds for all $c>0$, it is possible to choose a subsequence
of $\varepsilon$ (denote the subsequence by $\varepsilon$ again) and
positive numbers $\beta(\varepsilon)$ which satisfy $ \lim
_{\varepsilon\downarrow0} \beta(\varepsilon) = + \infty$, and
%
\begin{equation}\label{XY}
\lim_{\varepsilon\downarrow0}  \bigl| P_O^\varepsilon\bigl(w\bigl(T^{\beta
(\varepsilon)\varepsilon}\bigr)\in\Omega_i ^\varepsilon\bigr) -
Q_O^\varepsilon\bigl(w\bigl(T^{\beta(\varepsilon)\varepsilon}\bigr)\in\Omega
_i^\varepsilon\bigr) \bigr| =0
\end{equation}
for $i=1,2,\ldots,N$.
Let $\delta(\varepsilon) := \varepsilon\beta(\varepsilon)$.
Then, $\delta(\varepsilon)$ satisfies the conditions in Lemma~\ref{lm-es2}.

Now we assume that $\sigma(O)=I_n$ where $I_n$ means the unit matrix.
This assumption enables us to determine the weights of the edges explicitly.
Let
\[
p_i:= \frac{c_i ^{n-1} \int_0^1 r^{n-2}e^{-u_i(r)}\,dr}{\sum_{i=1}^N
c_i ^{n-1} \int_0^1 r^{n-2} e^{-u_i(r)}\,dr}.\vadjust{\goodbreak}
\]
We remark that when $u_i$ is independent of $i$, then we have\vspace*{-1pt}
$ p_i:= c_i ^{n-1} /  \break( \sum_{i=1}^N c_i ^{n-1}  ) $; hence
the weights $\{ p_i\}$ are determined by the ratio of the area of the
cross-section around the edge $I_i$.
Then, the following lemma holds.

\begin{lm}\label{lm-bc}
If $\sigma(O)=I_n$, then
\[
\lim_{\varepsilon\downarrow0} \sup_{|x|\leq3\kappa\varepsilon
} \bigl| P_x^\varepsilon\bigl(w\bigl(T^{\delta(\varepsilon)}\bigr)\in\Omega
_i^\varepsilon\bigr) -p_i \bigr| =0
\]
for $i=1,\ldots,N$.
\end{lm}

\begin{pf}
Applying Lemma~\ref{lm-es2} to both $X_{\cdot}^\varepsilon$ and
$Y_{\cdot}^\varepsilon$, and using~(\ref{XY}), it is sufficient to show that
%
\begin{equation}\label{lm-bc-7}
\lim_{\varepsilon\downarrow0}  \bigl| Q_O^\varepsilon\bigl(w\bigl(T^{\delta
(\varepsilon)}\bigr)\in\Omega_i^\varepsilon\bigr) -p_i \bigr| =0
\end{equation}
for $i=1,\ldots,N$.

We make a similar discussion as in the proof of Theorem 6.1 in~\cite{FW}.
Let~$\nu^{\varepsilon}$ be the invariant measure of the Markov chain
$\{ Y^\varepsilon(\tau^\varepsilon_{k})\}$, where $\tau^\varepsilon
_{k}$ are stopping times defined by
\begin{eqnarray*}
\tau_0 ^\varepsilon(w)&:=& 0,\\
\tilde\tau_{k} ^\varepsilon(w)&:=& \inf\{ u> \tau_{k-1}^\varepsilon
(w) \dvtx  |\pi(w(u))|> \delta(\varepsilon)\},\qquad k\in{\mathbb N},\\
\tau_{k} ^\varepsilon(w)&:=& \inf\{ u> \tilde\tau_k ^\varepsilon
(w) \dvtx  |\pi(w(u))|< 3\kappa\varepsilon\},\qquad k\in{\mathbb N}.
\end{eqnarray*}
Define a measure $\mu^\varepsilon$ on $\Omega^\varepsilon$ by
\[
\mu^\varepsilon(dx) := \exp(-U^\varepsilon(x))\,dx,\qquad x \in\Omega
^\varepsilon,
\]
a function space ${\mathscr D}({\mathscr E}^\varepsilon)$ by $ \{
f\in C^2 (\Omega^\varepsilon) \dvtx  \lim_{x\dvtx  d(x,\partial\Omega
^\varepsilon) \rightarrow0 }f(x)=0 \}$ and a~bilinear form $\mathscr
E ^\varepsilon$ by
\[
{\mathscr E}^\varepsilon(f,g) := \int_{\Omega^\varepsilon} \langle
\nabla f (x), \nabla g (x)\rangle\mu^\varepsilon(dx), \qquad f,g \in
{\mathscr D}({\mathscr E}^\varepsilon).
\]
Then, the pre-Dirichlet form $({\mathscr E}^\varepsilon, {\mathscr
D}({\mathscr E}^\varepsilon))$ on $L^2(\Omega^\varepsilon, \mu
^\varepsilon)$ is closable, and $Y^\varepsilon$ is associated to the
Dirichlet form obtained by closing $({\mathscr E}^\varepsilon,
{\mathscr D}({\mathscr E}^\varepsilon))$.
Note that $\mu^\varepsilon$ is an invariant measure of $Y^\varepsilon
$; see~\cite{FOT}.
By Theorem 2.1 in~\cite{Kh} we have
\[
\mu^\varepsilon(B) = \int_{\{x \in\Omega^\varepsilon\dvtx  |\pi(x)| =
3\kappa\varepsilon\}} \nu^{\varepsilon}(dx) \int \biggl[ \int
_0^{\tau_1^{\varepsilon}(w)} {\mathbb I}_B (w(t))\,dt  \biggr] Q_x^\varepsilon(dw)
\]
for $B \in{\mathscr B}({\mathbb R}^n)$.
Let $ B_i^{\varepsilon}:= \{ x \in\Omega_i^\varepsilon\dvtx  \delta
(\varepsilon)\leq|\pi(x)|\leq2\delta(\varepsilon)\}$.
Then,
%
\begin{eqnarray}\label{lm-bc-1}
&&   \mu^\varepsilon(B_i^{\varepsilon}) \nonumber\hspace*{-35pt}\\
&&   \qquad = \int_{\{x \in\Omega^\varepsilon\dvtx  |\pi(x)| = 3\kappa
\varepsilon\}} \nu^{\varepsilon} (dx) \int \!\biggl[ {\mathbb I}_{\Omega
_i} (w(\tilde{\tau}^{\varepsilon}_1)) \!\int_{\tilde\tau_1^{\varepsilon}(w)}^{\tau_1^{\varepsilon}(w)} {\mathbb
I}_{B_i^{\varepsilon}} (w(t))\,dt  \biggr] Q_x^\varepsilon(dw)\hspace*{-35pt}
\\
&&   \qquad = \int_{\{x \in\Omega^\varepsilon\dvtx  |\pi(x)| = 3\kappa
\varepsilon\}} \nu^{\varepsilon} (dx)  \int {\mathbb I}_{\Omega _i}
\bigl(w\bigl(T^{\delta (\varepsilon )}\bigr)\bigr) Q_x^\varepsilon (dw) \nonumber\hspace*{-35pt}\\
&& \hphantom{= \int_{\{x \in\Omega^\varepsilon\dvtx  |\pi(x)| = 3\kappa
\varepsilon\}} \nu^{\varepsilon} (dx) \int}
    \quad \times\int \biggl[ \int_{0}^{T^{3\kappa\varepsilon
}} {\mathbb I}_{B_i^{\varepsilon}} (\tilde w(t))\,dt  \biggr]
Q_{w(T^{\delta(\varepsilon)})}^\varepsilon(d\tilde w) .
\nonumber\hspace*{-35pt}
\end{eqnarray}
On the other hand, let
\[
Z(t):= -\delta(\varepsilon) + \check W(t), \qquad\check T := \inf\{
t>0\dvtx  |Z(t)|>2\delta(\varepsilon)-3 \kappa\varepsilon\},
\]
where $\check W$ is a one-dimensional Wiener process starting from 0, and
\[
F(x):= \int_{-2\delta(\varepsilon)}^x \int_{-2\delta(\varepsilon
)}^y {\mathbb I}_{[-\delta(\varepsilon),\delta(\varepsilon)]}(z)\,dz
\,dy, \qquad x \in{\mathbb R}.
\]
Then, by It\^o's formula we have
\[
E [ F(Z(\check T)) ] - F(-\delta(\varepsilon)) = \frac12
E \biggl[ \int_0^{\check T} {\mathbb I}_{[-\delta(\varepsilon),\delta
(\varepsilon)]}(Z_t)\,dt  \biggr] .
\]
Since $F$ can be computed explicitly, we see that $F(-\delta
(\varepsilon))=0$ and
\begin{eqnarray*}
E [ F(Z(\check T)) ] &=& F\bigl(2\delta(\varepsilon)-3 \kappa
\varepsilon\bigr) P\bigl ( Z(\check T)= 2\delta(\varepsilon)-3 \kappa
\varepsilon \bigr) \\
&=& \frac{\delta(\varepsilon)-3\kappa\varepsilon}{4\delta
(\varepsilon) -6\kappa\varepsilon} \bigl[ 2\delta(\varepsilon) ^2+
2\delta(\varepsilon) \bigl(\delta(\varepsilon) -3\kappa\varepsilon
\bigr) \bigr] .
\end{eqnarray*}
Thus, it follows that
\[
E \biggl[ \int_0^{\check T} {\mathbb I}_{[-\delta(\varepsilon),\delta
(\varepsilon)]}(Z(t))\,dt  \biggr] = 2\delta(\varepsilon) ^2 + o(\delta
(\varepsilon)^2).
\]
On the other hand, the strong Markov property and the reflection
principle imply that
\[
\int \biggl( \int_{0}^{T^{3\kappa\varepsilon}} {\mathbb
I}_{B_i^\varepsilon} (w(t))\,dt  \biggr) Q_y^\varepsilon(dw) = E \biggl[
\int_0^{\check T} {\mathbb I}_{[-\delta(\varepsilon),\delta
(\varepsilon)]}(Z(t))\,dt  \biggr]
\]
for all $y \in\{ x\in\Omega^\varepsilon\dvtx  |\pi(x)|=\delta
(\varepsilon)\}$, because the left-hand side is independent of the
behavior of $w$ moving in $\{ x\in\Omega^\varepsilon\dvtx  |\pi(x)|\geq
\delta(\varepsilon)\}$ under $Q_y^\varepsilon$.
Hence, it holds that
%
\begin{equation}\label{lm-bc-10}
\int \biggl( \int_{0}^{T^{3\kappa\varepsilon}} {\mathbb
I}_{B_i^\varepsilon} (w(t))\,dt  \biggr) Q_y^\varepsilon(dw)= 2 \delta
(\varepsilon)^2 + o(\delta(\varepsilon)^2 )
\end{equation}
for all $y \in\{ x\in\Omega^\varepsilon\dvtx  |\pi(x)|=\delta
(\varepsilon)\}$.
By Lemma~\ref{lm-es2},~(\ref{lm-bc-1}) and~(\ref{lm-bc-10}), we have
%
\begin{eqnarray}\label{lm-bc-2}
\mu^\varepsilon(B_i^\varepsilon) &=&  \bigl( 2 \delta(\varepsilon
)^2 + o(\delta(\varepsilon)^2 ) \bigr) \nu^{\varepsilon}\bigl (\{x \in
\Omega^\varepsilon\dvtx  |\pi(x)| = 3\kappa\varepsilon\}\bigr) \nonumber
\\[-8pt]
\\[-8pt]
&& {} \times \bigl( Q_O^\varepsilon\bigl(w\bigl(T^{\delta(\varepsilon
)}\bigr)\in\Omega_i^\varepsilon\bigr) +o_\varepsilon(1)  \bigr) .
\nonumber
\end{eqnarray}
Since $ \sum_{i=1} ^N Q_O^\varepsilon(w(T^{\delta(\varepsilon
)})\in\Omega_i^\varepsilon)=1$, we have, as $\varepsilon\downarrow0$
%
\begin{equation}\label{lm-bc-3}
\nu^{\varepsilon}\bigl (\{x \in\Omega^\varepsilon\dvtx  |\pi(x)| = 3\kappa
\varepsilon\}\bigr) = \frac12 \delta(\varepsilon)^{-2} \sum_{i=1} ^N
\mu^\varepsilon(B_i^\varepsilon) +o_\varepsilon(1).\vadjust{\goodbreak}
\end{equation}
Dividing both sides of~(\ref{lm-bc-2}) by those of~(\ref{lm-bc-3}), we
obtain that
\[
Q_O^\varepsilon\bigl(w\bigl(T^{\delta(\varepsilon)}\bigr)\in\Omega_i^\varepsilon
\bigr) = \frac{\mu^\varepsilon(B_i^\varepsilon)}{\sum_{i=1} ^N \mu
^\varepsilon(B_i^\varepsilon)} +o_\varepsilon(1).
\]
By the definition of $\mu^\varepsilon$, the continuity of $\sigma$
and $b$, and $\sigma(O)=I_n$, $\mu^\varepsilon(B_i^\varepsilon)$
can be expressed explicitly as
\[
\mu^\varepsilon(B_i^\varepsilon) = \omega_{n-2}\delta(\varepsilon
)c_i ^{n-1}\varepsilon^{n-1} \int_0^1 r^{n-2}e^{-u_i(r)}\,dr,
\]
where $\omega_{n-2}$ is the area of the $(n-2)$-dimensional unit sphere.
Therefore,~(\ref{lm-bc-7}) is proved.
\end{pf}

The statement in Lemma~\ref{lm-bc} can be improved as follows.

\begin{lm}\label{lm-bc2}
%
\[
\lim_{\delta' \downarrow0} \lim_{\varepsilon\downarrow0} \sup
_{|x|\leq3\kappa\varepsilon} \bigl| P_x^\varepsilon\bigl(w(T^{\delta
'})\in\Omega_i^\varepsilon\bigr) -p_i \bigr| =0
\]
for $i=1,\ldots,N$.
\end{lm}

\begin{pf}
In view of Lemma~\ref{lm-es2}, it is sufficient to show
\[
\lim_{\delta' \downarrow0} \lim_{\varepsilon\downarrow0} \bigl|
P_O^\varepsilon\bigl(w(T^{\delta'})\in\Omega_i^\varepsilon\bigr) - p_i
\bigr| = 0
\]
for $i=1,2,\ldots,N$.
Define stopping times $\{ \tau^{\varepsilon}_k , \tilde\tau^{\varepsilon}_k \}$ by
\begin{eqnarray*}
\tau_0 ^\varepsilon(w)&:=& 0,\\[-2pt]
\tilde\tau_{k} ^\varepsilon(w)&:=& \inf\{ u> \tau_{k-1}^\varepsilon
(w) \dvtx  |\pi(w(u))|> \delta(\varepsilon)\},\qquad k\in{\mathbb N} ,\\[-2pt]
\tau_{k} ^\varepsilon(w)&:=& \inf\{ u> \tilde\tau_k ^\varepsilon
(w) \dvtx  |\pi(w(u))|< 3\kappa\varepsilon\},\qquad k\in{\mathbb N} .
\end{eqnarray*}
By the strong Markov property, we have
%
\begin{eqnarray}\label
{lm-bc2-1}
&& P_O^\varepsilon\bigl(w(T^{\delta'})\in\Omega_i^\varepsilon\bigr) \nonumber\\[-2pt]
&& \qquad  = \sum_{k=1}^\infty\int{\mathbb I}_{\{ \tau_{k-1}^\varepsilon<
T^{\delta'}\}}(w)P_O^\varepsilon(dw)\\[-2pt]
&& \hphantom{\sum_{k=1}^\infty\int}\qquad  \quad {}  \times\int P_y ^\varepsilon(T^{\delta'}< T^{3 \kappa
\varepsilon}) {\mathbb I}_{\Omega_i^\varepsilon} (y) P _{w(\tau
_{k-1}^\varepsilon)}\bigl(w\bigl(T^{\delta(\varepsilon)}\bigr)\in dy\bigr) \nonumber
\end{eqnarray}
and
%
\begin{eqnarray} \label{lm-bc2-2}
p_i
&=& p_i \sum_{k=1}^\infty P_O^\varepsilon(\tau_{k-1}^\varepsilon
<T^{\delta'} <\tau_k^\varepsilon)\nonumber\\[-2pt]
&=& \sum_{k=1}^\infty\int{\mathbb I}_{\{ \tau_{k-1}^\varepsilon<
T^{\delta'}\}}(w)P_O^\varepsilon(dw) \\[-2pt]
&&\hphantom{\sum_{k=1}^\infty\int} {} \times\int p_i P_y ^\varepsilon(T^{\delta'}< T^{3 \kappa
\varepsilon}) P
_{w(\tau_{k-1}^\varepsilon)}\bigl(w\bigl(T^{\delta(\varepsilon)}\bigr)\in dy\bigr) \nonumber\vadjust{\goodbreak}
\end{eqnarray}
for $i=1,2,\ldots,N$.
Let $h_-^\varepsilon$ and $h_+^\varepsilon$ be functions on
$[0,\infty)$ given by
\begin{eqnarray*}
h_-^\varepsilon(z) &:=& \max_{i} \sup_{x \in \Omega _i^\varepsilon : |\pi (x)|=\max \{ z,3\kappa \varepsilon\} } \frac{2\langle e_i,b(x)\rangle}{|\sigma(x) ^T e_i|^2}, \\
h_+^\varepsilon(z)&:=& \min_{i} \inf_{x \in \Omega _i^\varepsilon : |\pi (x)|=\max \{ z,3\kappa \varepsilon\}}\frac{2\langle e_i,b(x)\rangle}{|\sigma(x) ^T e_i|^2},
\end{eqnarray*}
respectively.
Define functions $s_- ^\varepsilon$ and $s_+ ^\varepsilon$ on
$[0,\infty)$ by
\begin{eqnarray*}
s_-^\varepsilon(z) := \int_0^z \exp \biggl( -\int_0 ^{z'}
h_-^\varepsilon(z'')\,dz'' \biggr)\,dz',\\
s_+^\varepsilon(z) := \int_0^z \exp \biggl( -\int_0 ^{z'}
h_+^\varepsilon(z'')\,dz'' \biggr)\,dz',
\end{eqnarray*}
respectively.
Then, for $ y\in\{ x\in\Omega_i^\varepsilon\dvtx  |\pi(x)|= \delta
(\varepsilon)\}$ we have
\begin{eqnarray*}
&& \int s_-^\varepsilon\bigl( \bigl|\pi\bigl(w(T^{\delta'}\wedge T^{3 \kappa
\varepsilon})\bigr)\bigr|\bigr) P_y^\varepsilon(dw) - s_-^\varepsilon(\delta
(\varepsilon)) \\
&& \qquad  = \int s_-^\varepsilon\bigl( \langle e_i, w(T^{\delta'}\wedge T^{3
\kappa\varepsilon})\rangle\bigr) P_y^\varepsilon(dw) - s_-^\varepsilon(
\langle e_i,y\rangle) \\
&& \qquad  = -\frac12 \int \biggl[ \int_0^{T^{\delta'}\wedge T^{3 \kappa
\varepsilon}} h_-^\varepsilon(w(s)) |\sigma(w(s)) ^T e_i|^2   \\
&&\hphantom{-\frac12 \int \biggl[ \int_0^{T^{\delta'}\wedge T^{3 \kappa
\varepsilon}}} \qquad  \quad {}    \times\exp\biggl ( -\int_0 ^{w(s)}
h_-^\varepsilon(z')\,dz' \biggr)  \,ds\biggr] P_y^\varepsilon(dw) \\
&& \qquad  \quad {} + \int\biggl [ \int_0^{T^{\delta'}\wedge T^{3 \kappa
\varepsilon}} \langle e_i, b(w(s))\rangle\exp \biggl( -\int_0 ^{w(s)}
h_-^\varepsilon(z')\,dz' \biggr)\,ds  \biggr] P_y^\varepsilon(dw) \\
&& \qquad \leq0.
\end{eqnarray*}
Hence, it holds that
\[
s_-^\varepsilon( \delta') P_y^\varepsilon ( T^{\delta'}<T^{3
\kappa\varepsilon} ) + s_-^\varepsilon( 3\kappa\varepsilon)
P_y^\varepsilon ( T^{\delta'}>T^{3 \kappa\varepsilon} )
\leq s_-^\varepsilon(\delta(\varepsilon))
\]
for $ y\in\{ x\in\Omega^\varepsilon\dvtx  |\pi(x)|= \delta
(\varepsilon)\}$.
Since
\[
P_y^\varepsilon ( T^{\delta'}<T^{3 \kappa\varepsilon} ) +
P_y^\varepsilon ( T^{\delta'}>T^{3 \kappa\varepsilon} ) =1,
\]
we have
%
\begin{equation}\label{lm-bc2-3}
P_y^\varepsilon ( T^{\delta'}<T^{3 \kappa\varepsilon} )
\leq\frac{s_-^\varepsilon(\delta(\varepsilon))-s_-^\varepsilon(
3\kappa\varepsilon)}{s_-^\varepsilon( \delta')-s_-^\varepsilon(
3\kappa\varepsilon)}
\end{equation}
for $ y\in\{ x\in\Omega^\varepsilon\dvtx  |\pi(x)|= \delta
(\varepsilon)\}$.
Similarly we have
%
\begin{equation}\label{lm-bc2-4}
P_y^\varepsilon ( T^{\delta'}<T^{3 \kappa\varepsilon} )
\geq\frac{s_+^\varepsilon(\delta(\varepsilon))-s_+^\varepsilon(
3\kappa\varepsilon)}{s_+^\varepsilon( \delta')-s_+^\varepsilon(
3\kappa\varepsilon)}
\end{equation}
for $ y\in\{ x\in\Omega^\varepsilon\dvtx  |\pi(x)|= \delta
(\varepsilon)\}$.
Let ${\mathcal N}_{T^{\delta'}}( X_O^\varepsilon)$ be the number of
transitions of $X_O^\varepsilon$ from the set $ \{ x\in\Omega
^\varepsilon\dvtx  |\pi(x)|=3\kappa\varepsilon\} $ to the set $ \{
x\in\Omega^\varepsilon\dvtx  |\pi(x)|=\delta(\varepsilon) \} $ during
the time interval $[0,T^{\delta'}(X_O^\varepsilon)]$.
By Lemma~\ref{lm-bc},~(\ref{lm-bc2-1}),~(\ref{lm-bc2-2}), (\ref
{lm-bc2-3}) and~(\ref{lm-bc2-4}), we have
\begin{eqnarray*}
&& P_O^\varepsilon\bigl(w(T^{\delta'})\in\Omega_i^\varepsilon\bigr) - p_i \\
&& \qquad  \leq\frac{s_-^\varepsilon(\delta(\varepsilon))-s_-^\varepsilon
( 3\kappa\varepsilon)}{s_-^\varepsilon( \delta')-s_-^\varepsilon(
3\kappa\varepsilon)} \\
&& \qquad  \quad {} \times\sum_{k=1}^\infty\int P _{w(\tau_{k-1}^\varepsilon
)}\bigl(w\bigl(T^{\delta(\varepsilon)}\bigr)\in\Omega_i^\varepsilon\bigr) {\mathbb
I}_{\{ \tau_{k-1}^\varepsilon< T^{\delta'}\}}(w)P_O^\varepsilon(dw)
\\
&& \qquad  \quad {}  - \frac{s_+^\varepsilon(\delta(\varepsilon
))-s_+^\varepsilon( 3\kappa\varepsilon)}{s_+^\varepsilon( \delta
')-s_+^\varepsilon( 3\kappa\varepsilon)} p_i \sum_{k=1}^\infty\int
{\mathbb I}_{\{ \tau_{k-1}^\varepsilon< T^{\delta'}\}
}(w)P_O^\varepsilon(dw) \\
&& \qquad   \leq \biggl( \frac{s_-^\varepsilon(\delta(\varepsilon
))-s_-^\varepsilon( 3\kappa\varepsilon)}{s_-^\varepsilon( \delta
')-s_-^\varepsilon( 3\kappa\varepsilon)} - \frac{s_+^\varepsilon
(\delta(\varepsilon))-s_+^\varepsilon( 3\kappa\varepsilon
)}{s_+^\varepsilon( \delta')-s_+^\varepsilon( 3\kappa\varepsilon
)} \biggr) p_i E [ {\mathcal N}_{T^{\delta'}}( X_O^\varepsilon
) ] \\
&& \qquad  \quad {}  + o_\varepsilon(1) \frac{s_-^\varepsilon(\delta
(\varepsilon))-s_-^\varepsilon( 3\kappa\varepsilon
)}{s_-^\varepsilon( \delta')-s_-^\varepsilon( 3\kappa\varepsilon)}
E [ {\mathcal N}_{T^{\delta'}}( X_O^\varepsilon) ].
\end{eqnarray*}
By the definitions of $s_-^\varepsilon$ and $s_+^\varepsilon$, we obtain
\[
\limsup_{\varepsilon\downarrow0} \delta(\varepsilon)^{-1}
\biggl(
\frac{s_-^\varepsilon(\delta(\varepsilon))-s_-^\varepsilon(
3\kappa\varepsilon)}{s_-^\varepsilon( \delta')-s_-^\varepsilon(
3\kappa\varepsilon)}- \frac{s_+^\varepsilon(\delta(\varepsilon
))-s_+^\varepsilon( 3\kappa\varepsilon)}{s_+^\varepsilon( \delta
')-s_+^\varepsilon( 3\kappa\varepsilon)}  \biggr) \leq C,
\]
where $C$ is a constant independent of $\delta'$, and for each $\delta ' >0$
\[
\frac{s_-^\varepsilon(\delta(\varepsilon))-s_-^\varepsilon(
3\kappa\varepsilon)}{s_-^\varepsilon( \delta')-s_-^\varepsilon(
3\kappa\varepsilon)} = O(\delta(\varepsilon) ).
\]
On the other hand, a similar discussion as in the proof of Lemma \ref
{lm-es1} implies
\[
E [ {\mathcal N}_{T^{\delta'}}( X_O^\varepsilon) ] =
\delta (\varepsilon)^{-1} o_{\delta '}(1).
\]
Therefore, we have
\[
\limsup_{\delta' \downarrow0} \limsup_{\varepsilon\downarrow
0}\bigl(P_O^\varepsilon\bigl(w(T^{\delta'})\in\Omega_i^\varepsilon\bigr) - p_i\bigr)
\leq0.
\]
Similarly we obtain
\[
\limsup_{\delta' \downarrow0} \limsup_{\varepsilon\downarrow
0}\bigl(p_i - P_O^\varepsilon\bigl(w(T^{\delta'})\in\Omega_i^\varepsilon\bigr) \bigr)
\leq0.
\]
These inequalities yield the conclusion.
\end{pf}

We need a little more improvement of Lemma~\ref{lm-bc2} as follows.

\begin{lm}\label{lem3.7}
\[
\lim _{\delta ' \downarrow 0} \lim _{\delta \downarrow 0} \lim _{\varepsilon \downarrow 0}
\sup _{|x|\leq \delta }\bigl| P_x^\varepsilon \bigl(w(T^{\delta '})\in \Omega _i^\varepsilon \bigr) -p_i\bigr| =0
\]
for $i=1,\ldots ,N$.
\end{lm}

\begin{pf}
In view of Lemma~\ref{lm-bc2} it is sufficient to show
\[
\lim _{\delta ' \downarrow 0} \lim _{\delta \downarrow 0} \lim _{\varepsilon \downarrow 0}\sup _{3\kappa \varepsilon \leq |x|\leq \delta }
\bigl| P_x^\varepsilon \bigl(w(T^{\delta '})\in \Omega _i^\varepsilon \bigr) -p_i\bigr| =0.
\]
By Lemma~\ref{lm-bc2} again,
\begin{eqnarray*}
&& \bigl| P_x^\varepsilon \bigl(w(T^{\delta '})\in \Omega _i^\varepsilon \bigr) -p_i\bigr| \\
&&\qquad = \biggl| \int _{\{ y\in \Omega ^\varepsilon \dvtx |\pi (y)|= 3\kappa \varepsilon \}} P_y^\varepsilon \bigl(\tilde w(T^{\delta '})
\in \Omega _i^\varepsilon \bigr) P_x^\varepsilon \bigl(w(T^{3\kappa \varepsilon})\in dy , T^{3\kappa \varepsilon} < T^{\delta '}\bigr) \\
&& \hspace*{196pt}{} + P_x^\varepsilon ( T^{3\kappa \varepsilon} > T^{\delta '} ){\mathbb I}_{\Omega _i^\varepsilon}(x) -p_i\biggr| \\
&&\qquad =\bigl | \bigl(p_i + o_{\varepsilon ,\delta '}(1) \bigr) P_x^\varepsilon ( T^{3\kappa \varepsilon} < T^{\delta '})
+ P_x^\varepsilon ( T^{3\kappa \varepsilon} > T^{\delta '} ){\mathbb I}_{\Omega _i^\varepsilon}(x) -p_i| \\
&&\qquad \leq p_i | P_x^\varepsilon ( T^{3\kappa \varepsilon} < T^{\delta '} ) -1\bigr|
 +P_x^\varepsilon ( T^{3\kappa \varepsilon} > T^{\delta '} ){\mathbb I}_{\Omega _i^\varepsilon}(x) + o_{\varepsilon ,\delta '}(1).
\end{eqnarray*}
Here, $o_{\varepsilon ,\delta '}(1)$ means a term which converges to $0$ as $\delta '\downarrow 0$ after letting $\varepsilon \downarrow 0$.
Hence, it is sufficient to show for $\delta '> 0$ and $i=1,2,\ldots ,N$
\begin{equation}\label{eq3.28}
\lim _{\delta \downarrow 0} \lim _{\varepsilon \downarrow 0}\inf _{x\in \Omega _i^\varepsilon : 3\kappa \varepsilon \leq |x|\leq \delta }
 P_x^\varepsilon ( T^{3\kappa \varepsilon} < T^{\delta '} ) =1.
\end{equation}
Let $T^O(w):= \inf \{ t \geq 0 \dvtx w(t) =O\}$ and fix $i$.
By Theorem~\ref{th-curve} the law of $(T^{3\kappa \varepsilon }(X_{x^\varepsilon }^\varepsilon) , T^{\delta '}(X_{x^\varepsilon }^\varepsilon ))$
converges to that of $(T^O (X_x ) , T^{\delta '}(X_x ))$ as $\varepsilon \downarrow 0$ for
$x^\varepsilon \in \{ y\in \Omega ^\varepsilon \dvtx 3\kappa \varepsilon \leq |y|\leq \delta \}$ such that $x^\varepsilon$ converges to $x \in I_i$,
where the process $X_x$ is determined by the following stochastic differential equation:
\begin{eqnarray}
X_x(t) = x + \int _0^t \langle e_i, \sigma (X_x(s)) \,dW(s) \rangle + \int _0^t \langle e_i,b(X_x(s))\rangle \,ds,\nonumber\\
\eqntext{t \in [0,T^O (X_x ) \wedge T^{\delta '}(X_x )].}
\end{eqnarray}
By using $I_i=\bigcap _{\varepsilon '>0} \bigcup _{\varepsilon <\varepsilon '} \Omega _i^\varepsilon$ and compactness
of $\{ y\in {\mathbb R} ^n \dvtx |y|\leq \delta \}$, we have
\begin{eqnarray*}
&&\lim _{\varepsilon \downarrow 0} \inf _{x\in \Omega _i^\varepsilon : 3\kappa \varepsilon \leq |x|\leq \delta }
P_x^\varepsilon ( T^{3\kappa \varepsilon} < T^{\delta '} ) \\
&&\qquad= \inf _{x\in I_i : 0 \leq |x|\leq \delta }  P\bigl( T^O(X_x) < T^{\delta '}(X_x) \bigr).
\end{eqnarray*}
Since $\sigma \sigma ^T$ is uniformly positive definite, we have
\[
\lim _{\delta \downarrow 0} \inf _{x\in I_i : 0 \leq |x|\leq \delta }  P\bigl( T^O(X_x) < T^{\delta '}(X_x) \bigr) =1.
\]
This proves (3.28).\vadjust{\goodbreak}
\end{pf}


The lemmas above determine the boundary condition for $X$ at $O$.
Now let us characterize $X$ by a generator of a process on $\Gamma$.
Let
\[
\partial_{e_i}f(x):=\lim_{s\rightarrow0}\frac1s  \bigl(
f(x+se_i)-f(x) \bigr)
\]
for any differentiable function $f$ on $I_i$ and $i=1,2,\ldots,N$.
Define a second-order differential operator ${\mathcal L}_i$ on $I_i$ by
%
\begin{equation}\label{defLi}
{\mathcal L}_i := \tfrac12 |\sigma^T(x)e_i|^2 \,\partial_{e_i}^2 +
\langle b(x),e_i\rangle\,\partial_{e_i}
\end{equation}
for $i=1,2,\ldots,N$.
Define the second-order differential operator ${\mathcal L}$ on
$C_0(\Gamma)$
\begin{eqnarray*}
{\mathscr D}({\mathcal L}) &:=&  \Biggl\{ f\in C_0(\Gamma
)\dvtx  f|_{I_i\setminus O} \in C_b^2(I_i\setminus O)    \mbox{ for all }
i=1,2,\ldots,N,  \\
&&  \hphantom{\Biggl\{}  \lim_{s\downarrow0}{\mathcal L}_i f(se_i)    \mbox{ has
a common value for } i=1,2,\ldots,N,  \\
&& \hspace*{122pt}\hphantom{\Biggl\{}   \sum_{i=1} ^N p_i  \Bigl( \lim_{s\downarrow
0}(\partial_{e_i}f)(se_i) \Bigr) =0 \Biggr\} ,\\
{\mathcal L}f(x)&:=& {\mathcal L}_i f(x),  \qquad  x\in
I_i\setminus O ,\\
{\mathcal L}f(O)&:=& \lim_{s\downarrow0} {\mathcal L}_i f(se_i).
\end{eqnarray*}
Note that ${\mathcal L}f(O)$ does not depend on the selection of
$i=1,2,\ldots,N$.
We call $\{ p_i\}$ the weights of the Kirchhoff boundary condition at
$O$, and call $\sum_{i=1} ^N p_i  ( \lim_{s\downarrow0}(\partial
_{e_i}f)(se_i) ) =0$ the weighted Kirchhoff boundary condition at~$O$.

\begin{Th}\label{th-spider}
Consider diffusion processes $X^\varepsilon$ defined by~(\ref{SDE2}).
Assume that $\sigma(O)=I_n$ and the law of $X^\varepsilon(0)$ converges
to a probability measure $\mu_0$ on $\Gamma$.
Then, $X^\varepsilon$ converges weakly on $C([0,+\infty);{\mathbb
R}^n)$ to the diffusion process $X$ as $\varepsilon\downarrow0$,
where $X$ is determined by the conditions that the law of~$X(0)$ is
equal to $\mu_0$ and
%
\begin{equation}\label{th-spider-1}
E \biggl[   f(X(t))-f(X(s))-\int_s^t {\mathcal L}f(X(u))\,du  \Big|
{\mathscr F}_s  \biggr] =0
\end{equation}
for $t\geq s \geq0$ and $f \in{\mathscr D}({\mathcal L})$, where
$({\mathscr F}_t)$ is the filtration generated by $X$.
Therefore, ${\mathcal L}$ is the generator of $X$.
\end{Th}

\begin{pf}
From Lemma~\ref{lm-tight} we have that $\{ X^\varepsilon\}$ is tight.
We are going to show that there is a unique limit point in this family.
Let $X$ be any limit point of $\{ X^\varepsilon\}$, and denote the
sequence converging to $X$ by $\{ X^\varepsilon\}$ again.
Since this martingale problem is well-posed (see~\cite{FW};~\cite{EK}
for the relationship between martingale problems and partial
differential equations, and~\cite{Lu} for the uniqueness of the
semigroup generated by $\mathcal{L}$), it is sufficient to prove that\vadjust{\goodbreak}
$X$ satisfies~(\ref{th-spider-1}).
Fix $s\geq0$. Let $\delta'$ be a positive number.
Define the following stopping times:
\begin{eqnarray*}
\tilde\tau_{0} &:=& s, \\[-2pt]
\tau_0 &:=& \inf\{ u\geq s \dvtx  X(u) =O \},\\[-2pt]
\tilde\tau_{k} &:=& \inf\{ u> \tau_{k-1} \dvtx  |X(u)|> \delta'\}, \qquad
k\in{\mathbb N},\\[-2pt]
\tau_{k} &:=& \inf\{ u> \tilde\tau_k \dvtx  X(u)=O \}, \qquad  k\in{\mathbb N}.
\end{eqnarray*}
Then, for $f \in{\mathscr D}({\mathcal L})$, $s \leq t$
\begin{eqnarray*}
&& E \biggl[   f(X(t))-f(X(s))-\int_s^t {\mathcal L}f(X(u))\,du
\Big| {\mathscr F}_s  \biggr] \\[-2pt]
&& \qquad  = E\Biggl [   \sum_{k=1}^{\infty}  \biggl( f\bigl(X(t\wedge\tilde\tau
_{k})\bigr)-f\bigl(X(t \wedge\tau_{k-1})\bigr)-\int_{t \wedge\tau_{k-1}}^{t\wedge
\tilde\tau_{k}} {\mathcal L}f(X(u))\,du  \biggr)  \Big| {\mathscr F}_s
 \Biggr] \\[-2pt]
&& \qquad  \quad {}  + \sum_{k=0}^{\infty} E\biggl [   f\bigl(X(t\wedge\tau
_{k})\bigr)-f\bigl(X(t \wedge\tilde\tau_{k})\bigr)-\int_{t \wedge\tilde\tau
_{k}}^{t\wedge\tau_{k}} {\mathcal L}f(X(u))\,du  \Big| {\mathscr F}_s
 \biggr] .
\end{eqnarray*}
Because of Theorem~\ref{th-curve} the second sum vanishes.
We estimate the first sum as follows:
\begin{eqnarray*}
&&\Biggl | E\Biggl[ \sum _{k=1}^{\infty} \biggl( f\bigl(X(t\wedge \tilde \tau _{k})\bigr)
-f\bigl(X(t \wedge \tau _{k-1})\bigr)-\int _{t \wedge \tau _{k-1}}^{t\wedge \tilde \tau _{k}}
{\mathcal L}f(X(u))\,du \biggr) \Biggr] \Big| {\mathscr F}_s\Biggr|\\[-2pt]
&&\qquad \leq
\biggl| E\biggl[ \sum _{k: \tilde \tau _k < t}
\bigl( f(X(\tilde \tau _{k}))-f(X(\tau _{k-1})) \bigr) \Big| {\mathscr F}_s\biggr]
\biggr|\\[-2pt]
&& \qquad  \quad {}  + \|{\mathcal L}f\|_{\infty} E \biggl[ \int_{s}^{t} {\mathbb
I}_{\{ x\dvtx  |x|\leq\delta'\}}(X(u))\,du  \Big| {\mathscr F}_s\biggr] + \sup _{|x|\leq \delta '}|
f(x)-f(O)|.
\end{eqnarray*}
Clearly, the third term on the right-hand side converges to $0$ as
$\delta' \downarrow0$.
By Lemma~\ref{lm-es1} the second term on the right-hand side converges
to $0$ as $\delta'\downarrow0$.
The first sum on the right-hand side is equal to
%
\begin{eqnarray} \label{th2-1-1}\qquad
&&\Biggl| \sum _{k=1}^{\infty} \sum _{i=1}^N \bigl( f(\delta 'e_i)-f(O)  \bigr)
 P\bigl(X(\tilde \tau _{k})\in I_i, \tilde \tau _k < t| {\mathscr F}_s\bigr)\biggr|
 \nonumber
 \\[-9pt]
 \\[-9pt]
 \nonumber
&&\qquad =\Biggl| \sum _{k=1}^{\infty} \sum _{i=1}^N \Bigl( \delta '\lim _{s\downarrow 0} f'(se_i)
+ o(\delta ')\Bigr) P\bigl(X(\tilde \tau _{k})\in I_i, \tilde \tau _k < t| {\mathscr F}_s\bigr) \biggr|
\end{eqnarray}
%
Let $\delta \in (0,\delta ')$ and let, for any $\varepsilon >0$:
\begin{eqnarray*}
\tau _0^{\varepsilon ,\delta }&:=& \inf\{ u> s \dvtx |\pi (X^\varepsilon (u))| < \delta \},\\[-2pt]
\tilde \tau _{k} ^{\varepsilon ,\delta }&:=& \inf\{ u> \tau _{k-1} ^{\varepsilon ,\delta }:
|\pi (X^\varepsilon (u))|> \delta '\},\qquad k\in {\mathbb N},\\[-2pt]
\tau _{k} ^{\varepsilon ,\delta }&:=& \inf\{ u> \tilde \tau _k ^{\varepsilon ,\delta }\dvtx
|\pi (X^\varepsilon (u))| < \delta \},\qquad k\in {\mathbb N}.
\end{eqnarray*}
The distributions of the pairs $( X^\varepsilon , \tilde \tau _{k} ^{\varepsilon ,\delta },
\tau _{k} ^{\varepsilon ,\delta })$ converge weakly to those of $( X, \tilde \tau _{k} , \tau _{k} )$ as $\delta \downarrow 0$ after $\varepsilon \downarrow 0$.
Hence, by Lemma~\ref{lem3.7} we have
\begin{eqnarray*}
&& P\bigl(X(\tilde \tau _{k})\in I_i, \tilde \tau _{k} < t| {\mathscr F}_s\bigr)\\[-2pt]
&&\qquad = \lim _{\delta \downarrow 0} \lim _{\varepsilon \downarrow 0} P\bigl(X^\varepsilon (\tilde \tau _{k} ^{\varepsilon ,\delta })
\in \Omega _i^\varepsilon , \tilde \tau _{k} ^{\varepsilon ,\delta } < t| {\mathscr F}_s\bigr) \\[-2pt]
&&\qquad = \lim _{\delta \downarrow 0} \lim _{\varepsilon \downarrow 0} \int P^\varepsilon _y \bigl(w(T^{\delta '})\in \Omega _i^\varepsilon \bigr)
P\bigl(X^\varepsilon (\tau ^{\varepsilon ,\delta } _{k-1})\in dy, \tilde \tau ^{\varepsilon ,\delta }_k < t| {\mathscr F}_s\bigr) \\[-2pt]
&&\qquad = \bigl( p_i + o_{\delta '}(1) \bigr) P( \tilde \tau _{k} < t| {\mathscr F}_s) .
\end{eqnarray*}
%
Note that $ \sum_{k=1}^{\infty} P( \tilde{\tau}_{k} < t| {\mathscr F}_s)$ is equal to
the expectation of the number~of transitions of $X$ from the point $O$
to the set $ \{ x\in\Gamma\dvtx  |x|= \delta' \} $ during the time
interval $[s,t]$ [with respect to a general initial condition $X(0)$].
Approximating that by the expectation of the number of transitions of
$X^\varepsilon$ from the set $ \{ x\in\Omega^\varepsilon\dvtx  |\pi
(x)|= \delta\} $ to the set $ \{ x\in\Omega
^\varepsilon\dvtx  |\pi(x)|= \delta' \} $ during the time interval
$[s,t]$, similarly as in the proof of Lemma~\ref{lm-es1} we obtain the estimate
\[
\sum_{k=1}^{\infty} P( \tilde{\tau}_{k} < t| {\mathscr F}_s) \leq\frac{C_t}{\delta'}
\]
with a positive constant $C_t$ depending only on $t$.
Hence, by~(\ref{th2-1-1}) we have
\begin{eqnarray*}
&&\biggl| E\biggl[ \sum_{k :\tilde \tau _{k} < t}
\bigl( f\bigl(X(t\wedge \tilde \tau _{k})\bigr)-f\bigl(X(t\wedge \tau _{k-1})\bigr) \bigr)
\Big| {\mathscr F}_s\biggr] \biggr|  \\[-2pt]
&&\qquad\leq \frac {C_t}{\delta '}
\Biggl| \sum _{i=1}^N \delta '\lim _{s\downarrow 0} f'(se_i) p_i +o(\delta ') \Biggr|.
\end{eqnarray*}
Since $f\in{\mathscr D}({\mathcal L})$, the right-hand side converges
to $0$ as $\delta' \downarrow0$.
\end{pf}

Similarly as in Section~\ref{section curve}, the argument above is also
available in the case where the boundary of $\Omega^\varepsilon$
carries a Neumann boundary condition.
Consider a diffusion process $X^\varepsilon$ which is associated to $L$
in $\Omega^\varepsilon$ and satisfies the reflecting boundary
condition on $\partial\Omega^\varepsilon$.
Then, $X^\varepsilon$ can be expressed by the following equation:
%
\begin{equation}\label{SDE2-2}
  \widehat X^\varepsilon(t)=\widehat X^\varepsilon(0) + \int_0
^{t}\sigma(\widehat X^\varepsilon(s))\,dW(s) + \int_0^{t}b(\widehat
X^\varepsilon(s))\,ds + \Phi^\varepsilon(\widehat X^\varepsilon)(t),\hspace*{-35pt}
\end{equation}
where $\Phi^\varepsilon$ is a singular drift which forces the process
to be reflecting on~$\partial\Omega^\varepsilon$; see~\cite{SV2}.
Note that $\widehat X^\varepsilon$ depends on $\Omega^\varepsilon$ but
is independent of $U^\varepsilon$.
Discussing this case in a similar way as we did in the case of
Dirichlet boundary condition we obtain the following theorem.
Let
\begin{eqnarray*}
\widehat p_i&:=& \frac{c_i ^{n-1}}{\sum_{i=1}^N c_i ^{n-1}},\\[-1pt]
{\mathscr D}(\widehat{\mathcal L}) &:=&  \Biggl\{ f\in
C_0(\Gamma)\dvtx  f|_{I_i\setminus O} \in C_b^2(I_i\setminus O)    \mbox{ for
all } i=1,2,\ldots,N,  \\[-1pt]
&& \hphantom{\Biggl\{}   \lim_{s\downarrow0}{\mathcal L}_i f(se_i)    \mbox{ has
a common value for } i=1,2,\ldots,N,  \\[-1pt]
&&\hspace*{122pt}  \hphantom{\Biggl\{}  \sum_{i=1} ^N {\widehat p_i} \Bigl ( \lim_{s\downarrow
0}(\partial_{e_i}f)(se_i) \Bigr) =0 \Biggr\} ,\\
\widehat{\mathcal L}f(x)&:=& {\mathcal L}_i f(x),  \qquad
x\in I_i\setminus O ,\\
\widehat{\mathcal L}f(O)&:=& \lim_{s\downarrow0}
{\mathcal L}_i f(se_i),
\end{eqnarray*}
 where ${\mathcal L}_i$ is given by~(\ref{defLi}).
Note that $\widehat{\mathcal L}f(O)$ does not depend on the selection
of $i=1,2,\ldots,N$.

\begin{Th}\label{th-spider2}
Consider the diffusion processes\vspace*{1pt} $\widehat X^\varepsilon$ defined by
(\ref{SDE2-2}).
Assume that $\sigma(O)=I_n$ and the law of $\widehat X^\varepsilon(0)$
converges to a probability measure $\mu_0$ on~$\Gamma$.
Then, $\{ \widehat X^\varepsilon\}$ converge weakly on $C([0,+\infty
);{\mathbb R}^n)$ to the diffusion process $\widehat X$ as $\varepsilon
\downarrow0$, where $\widehat X$ is determined by the conditions that
the law of~$\widehat X(0)$ is equal to $\mu_0$ and
\[
E \biggl[   f(\widehat X(t))-f(\widehat X(s))-\int_s^t \widehat
{\mathcal L}f(\widehat X(u))\,du \Big | {\mathscr F}_s  \biggr] =0
\]
for $t\geq s \geq0$ and $f \in{\mathscr D}(\widehat{\mathcal L})$,
where $({\mathscr F}_t)$ is the filtration generated by $\widehat X$.
Therefore, $\widehat{\mathcal L}$ is the generator of $\widehat X$.
\end{Th}

\begin{re}
The weights $\{ \widehat p_i\}$ of the case of Neumann boundary
condition can be obtained from the wights $\{ p_i\}$ discussed in
Theorem~\ref{th-spider} in the heuristic limit where the potential
$u_i$ around each edge takes only the value~$0$ on $[0,1)$ and $+\infty
$ on $[1,+\infty)$.
\end{re}

\begin{re}\label{re-spider}
As mentioned in Remark~\ref{re-curve}, we can discuss similarly the
case where the shapes of the tubes $\{ \Omega_i^\varepsilon\}$ are
not cylindrical.
However, if $U^\varepsilon$ is not defined by a scaling of a fixed
function $U$, the weights of the weighted Kirchhoff boundary condition
cannot be determined uniquely.
To handle this more general case, we have to assume that $U^\varepsilon
$ satisfies some uniform bound.
\end{re}

\section{The case of general graphs}\label{section graph}

In this section we present results obtained by combining the results of
Sections~\ref{section curve} and~\ref{section spider}, and, in this
way, we cover more general graphs.
Let $\Lambda$ be a finite or countable set, $\Xi$ be a subset of
$\Lambda\times\Lambda$, $\{V_\lambda\dvtx  \lambda\in\Lambda\}$ be
vertices in ${\mathbb R}^n$, $\{ E_{\lambda,\lambda'}\dvtx  (\lambda
,\lambda') \in\Xi\}$ be $C^3$-curves with ends $\{ V_{\lambda},
V_{\lambda'}\}$ and $ G := \bigcup_{(\lambda,\lambda')\in\Xi}
E_{\lambda,\lambda'}$.
Denote $\lambda\sim\lambda'$ if $(\lambda,\lambda') \in\Xi$.\vadjust{\goodbreak}

Let us denote the length of $E_{\lambda,\lambda'}$ by $|E_{\lambda
,\lambda'}|$.
Define $(\gamma_{\lambda, \lambda'}(s)\dvtx s\in[0,|E_{\lambda, \lambda
'}|])$ as the arc-length parameterization of $E_{\lambda, \lambda'}$
with $\gamma_{\lambda, \lambda'}(0)=V_\lambda$.
Assume that the number of $\{V_\lambda\dvtx  \lambda\in\Lambda\} \cap\{x
\in{\mathbb R}^n\dvtx  |x|\leq M\}$ is finite for all $M>0$, $|E_{\lambda
,\lambda'}|$ is finite for all $(\lambda,\lambda') \in\Xi$ and
\[
\lim_{s\downarrow0}\langle\dot\gamma_{\lambda, \lambda_1}(s) ,
\dot\gamma_{\lambda, \lambda_2}(s)\rangle<1
\]
for all $\lambda\sim\lambda_1$ and $\lambda\sim\lambda_2$ such
that $\lambda_1 \neq\lambda_2$.
Let $c_{\lambda,\lambda'}$ be a positive number for $(\lambda
,\lambda') \in\Xi$, and let
\begin{eqnarray*}
\kappa_\lambda&:=& \max\Bigl \{ 2\sqrt2 c_{\lambda,\lambda_1}\big/\sqrt
{1-\lim_{s\downarrow0}\langle\dot\gamma_{\lambda, \lambda_1}(s)
, \dot\gamma_{\lambda, \lambda_2}(s)\rangle} \dvtx
\lambda_1 ,\lambda_2 \in\Lambda \\
&&\hspace*{79pt}\hphantom{\max\Bigl \{}  \mbox{such that } \lambda\sim\lambda_1 , \lambda\sim\lambda_2, \lambda _1 \neq \lambda _2  \Bigr\}
\end{eqnarray*}
for $\lambda\in\Lambda$.
Let $\pi(x)$ be a point in $G$ which is nearest to $x \in{\mathbb R}^n$.
Assume that there exists a small $\varepsilon_0 >0$ and positive numbers $\{ \kappa _\lambda\}$ such that $\pi
(x)$ is uniquely determined for all $ x \in\bigcup_{\lambda\sim
\lambda'} \{ x\in{\mathbb R}^n\dvtx d(x,E_{\lambda ,\lambda '}) <c_{\lambda ,\lambda '}\varepsilon, d(x,V_{\lambda}) \geq \kappa _\lambda \varepsilon
 \mbox{ and }
d(x,V_{\lambda '}) \geq \kappa _{\lambda '}\varepsilon \}$ and for all
$\varepsilon\in(0,\varepsilon_0]$, and that $\ddot\gamma_{\lambda
, \lambda'} (s)=0$ for sufficiently small $s$ for each $(\lambda
,\lambda') \in\Xi$.

Let $u_{\lambda,\lambda'}$ be given similarly to $u$ in Section \ref
{section curve} for $(\lambda,\lambda') \in\Xi$.
For $\varepsilon\in(0,\varepsilon_0]$, let $U^\varepsilon$ be a
function on ${\mathbb R}^n$ with values in $[0,+\infty]$, and assume
\begin{eqnarray}
U^\varepsilon(x) &= & u_{\lambda,\lambda'}(c_{\lambda,\lambda'}
^{-1}\varepsilon^{-1}d(x, E_{\lambda,\lambda'})),\nonumber\\
 \eqntext{x \in\{ x\in{\mathbb R}^n\dvtx  \pi(x)\in E_{\lambda,\lambda'},
d(x,E_{\lambda,\lambda'}) <c_{\lambda,\lambda'}\varepsilon,
d(x,V_{\lambda}) \geq\kappa_\lambda\varepsilon, d(x,V_{\lambda'}) \geq \kappa _{\lambda'} \varepsilon\} ,}\\
U^\varepsilon(x) &= & +\infty,\nonumber\\
 \eqntext {x \in\{ x\in{\mathbb R}^n\dvtx  \pi(x)\in E_{\lambda,\lambda'},
d(x,E_{\lambda,\lambda'}) \geq c_{\lambda,\lambda'}\varepsilon,
d(x,V_{\lambda}) \geq\kappa_\lambda\varepsilon, d(x,V_{\lambda'}) \geq \kappa _{\lambda'} \varepsilon\} ,}
\end{eqnarray}
$\Omega^\varepsilon:= \{ x\dvtx  U^\varepsilon(x)<\infty\}$ is a simply
connected domain, $ \partial\Omega^\varepsilon$ is an
$(n-1)$-dimensional $C^2$-manifold embedded in ${\mathbb R}^n$ and
$U^\varepsilon|_{\Omega^\varepsilon}$ is a $C^1$-function on~$\Omega
^\varepsilon$.
In addition, we assume
\[
\lim_{m\rightarrow \infty} \langle -\nabla U(x_m), \nabla d(x_m , \partial \Omega ^\varepsilon )\rangle
= +\infty \quad \mbox{and}\quad -\lim _{m\rightarrow \infty}\frac {U^\varepsilon (x_m)}{\operatorname{log}(d(x_m,\partial \Omega ^\varepsilon ))}
 = + \infty
\]
for any sequence $\{ x_m\}$ which converges to a point $x \in\partial
\Omega^\varepsilon$.

Consider a diffusion process $X^\varepsilon$ given by the following equation:
%
\begin{eqnarray}\label{SDE3}
X^\varepsilon(t)&=&X^\varepsilon(0) + \int_0 ^t \sigma(X^\varepsilon
(s))\,dW(s) + \int_0^t b(X^\varepsilon(s))\,ds \nonumber
\\[-8pt]
\\[-8pt]
&&   {}- \int_0^t (\nabla
U^\varepsilon)(X^\varepsilon(s))\,ds,
\nonumber
\end{eqnarray}
where $X^\varepsilon(0)$ is an $\Omega^\varepsilon$-valued random
variable, $W$ is an $n$-dimensional Wiener process, $ \sigma\in
C_b({\mathbb R}^n; {\mathbb R}^n\otimes{\mathbb R}^n) $ and $ b \in
C_b({\mathbb R}^n; {\mathbb R}^n)$.
Let $a:=\sigma\sigma^T$, and assume that $a$\vadjust{\goodbreak} is uniformly positive definite.
Define a second-order elliptic differential operator $L$ on $\Omega
^\varepsilon$ by
\[
L:= \frac12 \sum_{i,j=1}^n a_{ij}(x)\,\frac{\partial}{\partial
x_i}\,\frac{\partial}{\partial x_j}+ \sum_{i=1}^n b_i(x)\,\frac
{\partial}{\partial x_i}.
\]
Then $X^\varepsilon$ is associated with $(L-\langle\nabla
U^\varepsilon,\nabla\rangle)$.
Similarly to Section~\ref{section spider}, it holds that $X^\varepsilon
$ does not exit from $\Omega^\varepsilon$ almost surely.
Assume that $\sigma(V_{\lambda})= \sigma_\lambda I_n$ for all
$\lambda\in\Lambda$ where $\sigma_\lambda >0$.

For $(\lambda,\lambda')\in\Xi$, define a second-order differential
operator $\mathcal L_{\lambda,\lambda'}$ on $E_{\lambda,\lambda'}$~by
\begin{eqnarray*}
&&{\mathcal L}_{\lambda ,\lambda '} f(x)\\[-2pt]
&&\qquad:= \frac 12 |\sigma (x)^T \dot \gamma _{\lambda ,\lambda '} \circ \gamma _{\lambda ,\lambda '} ^{-1}(x)|^2
\frac {d^2}{ds^2}(f\circ \gamma _{\lambda ,\lambda '})(\gamma _{\lambda ,\lambda '}^{-1}(x)) \\[-2pt]
&&\qquad\quad{}+ [ \langle b(x),
\dot \gamma _{\lambda ,\lambda '}\circ \gamma _{\lambda ,\lambda '} ^{-1}(x)\rangle \\[-2pt]
&&\hspace*{27pt}{}\qquad  + \langle \sigma (x)^T \ddot \gamma _{\lambda ,\lambda '}
 \circ \gamma _{\lambda ,\lambda '} ^{-1}(x),\\[-2pt]
&&\hspace*{48pt}\qquad\sigma (x)^T \dot \gamma _{\lambda ,\lambda '} \circ \gamma _{\lambda ,\lambda '} ^{-1}(x)\rangle ]
\frac {d}{ds} (f\circ \gamma _{\lambda ,\lambda '})(\gamma _{\lambda ,\lambda '}^{-1}(x)),
\end{eqnarray*}
for $x\in E_{\lambda ,\lambda '}$ and $f\in C^2_b ( E_{\lambda ,\lambda '})$ where $s$ is the parameter for the arc-length
parametrization $\gamma _{\lambda ,\lambda '}$.
Let
\[
p_{\lambda,\lambda'} := \frac{ c_{\lambda,\lambda'}
^{n-1} \int_0^1 r^{n-2}\exp ( -u_{\lambda,\tilde\lambda
}(r) )\,dr}{\sum_{\tilde\lambda\dvtx \tilde\lambda
\sim\lambda}c_{\lambda,\tilde\lambda} ^{n-2} \int_0^1 r^{n-1}
\exp ( -u_{\lambda,\tilde\lambda}(r) )\,dr}.
\]
By using these notations, define the second-order differential operator
$\mathcal L$ on~$C_0(G)$ by
\begin{eqnarray*}
{\mathscr D}({\mathcal L}) &:=&  \biggl\{ f\in
C_0(G)\dvtx    f|_{E_{\lambda,\lambda'} \setminus\{ V_\lambda,V_{\lambda
'} \} } \in C_b^2(E_{\lambda,\lambda'} \setminus\{ V_\lambda
,V_{\lambda'} \} )    \mbox{ for } \lambda\sim\lambda',\\[-2pt]
&& \hphantom{\biggl\{} \mbox{for } \lambda\in\Lambda, \lim_{s\downarrow
0}{\mathcal L}_{\lambda,\lambda'} f(\gamma_{\lambda,\lambda'}(s))
   \mbox{ has a common value  for } \lambda' \dvtx  \lambda\sim\lambda',\\[-2pt]
&&\hspace*{78pt}    \sum_{\lambda'\dvtx \lambda'\sim\lambda} p_{\lambda
,\lambda'} \lim_{s\downarrow0} \biggl( \frac{d}{ds} \bigl(f\circ\gamma
_{\lambda,\lambda'} (s)\bigr) \biggr) =0 \mbox{ for } \lambda\in\Lambda
 \biggr\} ,\\[-2pt]
{\mathcal L}f(x)&:=& {\mathcal L}_{\lambda,\lambda'}
f(x),  \qquad  x\in E_{\lambda,\lambda'}, (\lambda,\lambda') \in\Xi
,\\[-2pt]
{\mathcal L}f(V_\lambda)&:=& \lim_{x \rightarrow
V_\lambda} {\mathcal L}_{\lambda,\lambda'} f(x),  \qquad \lambda\in
\Lambda,
\end{eqnarray*}
where the limit $x \rightarrow V_\lambda$ is along $E_{\lambda,\lambda'}$.
Note that ${\mathcal L}f(V_\lambda)$ does not depend on the selection
of $\lambda'$.

Since by locality the behavior of diffusion processes associated with
differential operators is determined in a given point by the behavior
in its neighborhoods, we have the following theorem by Theorem \ref
{th-curve} and~\ref{th-spider}.

\begin{Th}\label{th-graph}
Consider the diffusion process $X^\varepsilon$ defined by~(\ref{SDE3}).
Assume that the law of $X^\varepsilon(0)$ converges\vadjust{\goodbreak} to a probability
measure $\mu_0$ on $G$.
Then, $\{ X^\varepsilon\}$ converge weakly on $C([0,+\infty);
{\mathbb R}^n)$ to the diffusion process $X$ as $\varepsilon\downarrow
0$, where $X$ determined by the conditions that the law of $X(0)$ is
equal to $\mu_0$ and
\[
E \biggl[   f(X(t))-f(X(s))-\int_s^t {\mathcal L}f(X(u))\,du \Big |
{\mathscr F}_s  \biggr] =0
\]
for $t\geq s \geq0$ and all $f \in{\mathscr D}({\mathcal L})$, where
$({\mathscr F}_t)$ is the filtration generated by $X$.
The operator ${\mathcal L}$ as defined above is thus the generator of $X$.
\end{Th}

Similarly as in Sections~\ref{section curve} and~\ref{section spider},
our discussion is also available for the case where the boundary
$\Omega^\varepsilon$ carries a Neumann boundary condition for the process.
Consider a diffusion process $\widehat X^\varepsilon$ which is
associated with $L$ in $\Omega^\varepsilon$ and reflecting on
$\partial\Omega^\varepsilon$ [defined similarly as the process
described by~(\ref{SDE2-2})].

Let
\begin{eqnarray*}
\widehat p_{\lambda,\lambda'} &:=& \frac{ c_{\lambda
,\lambda'} ^{n-1} }{\sum_{\tilde\lambda\dvtx \tilde
\lambda\sim\lambda}c_{\lambda,\tilde\lambda} ^{n-1}},\\
{\mathscr D}(\widehat{\mathcal L}) &:=&  \biggl\{ f\in
C_0(G)\dvtx    f|_{E_{\lambda,\lambda'} \setminus\{ V_\lambda,V_{\lambda
'} \} } \in C_b^2(E_{\lambda,\lambda'} \setminus\{ V_\lambda
,V_{\lambda'} \} ) \mbox{ for } \lambda\sim\lambda',\\
&&\hphantom{\biggl\{}  \mbox{for } \lambda\in\Lambda, \lim_{s\downarrow
0}{\mathcal L}_{\lambda,\lambda'} f(\gamma_{\lambda,\lambda'}(s))
\mbox{ has a common value for } \lambda' \dvtx  \lambda\sim\lambda',\\
&&\hspace*{78pt}    \sum_{\lambda'\dvtx \lambda'\sim\lambda} \widehat
p_{\lambda,\lambda'} \lim_{s\downarrow0} \biggl( \frac{d}{ds} (f\circ
\gamma_{\lambda,\lambda'} (s)) \biggr) =0 \mbox{ for } \lambda\in
\Lambda \biggr\} ,\\
\widehat{\mathcal L}f(x)&:=& {\mathcal L}_{\lambda
,\lambda'} f(x),  \qquad  x\in E_{\lambda,\lambda'}, (\lambda,\lambda
') \in\Xi,\\
\widehat{\mathcal L}f(V_\lambda)&:=& \lim_{x
\rightarrow V_\lambda} {\mathcal L}_{\lambda,\lambda'} f(x),  \qquad
\lambda\in\Lambda,
\end{eqnarray*}
where the limit $x \rightarrow V_\lambda$ is along $E_{\lambda,\lambda'}$.
Then, we obtain the following theorem.

\begin{Th}\label{th-graph2}
Consider the diffusion process $\widehat X^\varepsilon$ defined above.
Assume that the law of $\widehat X^\varepsilon(0)$ converges to $\mu_0$.
Then, $\{ \widehat X^\varepsilon\}$ converge weakly on\vspace*{1pt} $C([0,+\infty
);{\mathbb R}^n)$ to the diffusion process $\widehat X$ as $\varepsilon
\downarrow0$, where $\widehat X$ is determined by the conditions\vspace*{1pt} that
the law of $\widehat X(0)$ is equal to $\mu_0$ and
\[
E \biggl[   f(\widehat X(t))-f(\widehat X(s))-\int_s^t \widehat
{\mathcal L}f(\widehat X(u))\,du \Big | {\mathscr F}_s  \biggr] =0
\]
for $t\geq s \geq0$ and $f \in{\mathscr D}(\widehat{\mathcal L})$
where $({\mathscr F}_t)$ is the filtration\vspace*{1pt} generated by $\widehat X$.
The operator $\widehat{\mathcal L}$ as defined above is thus the
generator of $\widehat X$.
\end{Th}

\begin{re}
As mentioned in Remarks~\ref{re-curve} and~\ref{re-spider}, similar
discussions can be given for the case where the shapes of the tubes are
not cylindrical.\vadjust{\goodbreak}
In the case where $\sigma=I_n$, $b=0$, and $E_{\lambda,\lambda'}$
are straight, the result of Theorem~\ref{th-graph2} coincides with
Theorem 6.1 in~\cite{FW}.
\end{re}

\section*{Acknowledgments}
We are very grateful to Claudio Cacciapuoti, Gianfausto Dell'Antonio,
Kazumasa Kuwada, Michael R\"ockner, Luciano Tubaro, Yohsuke Imagi and Hirokazu Maruhashi for very
interesting and stimulating discussions.

The first author would like to express his gratitude to Gianfausto
Dell'An\-tonio for his warm hospitality at SISSA (Trieste).
We also thank Luciano Tubaro, Raul Serapioni and Luca di Persio,
respectively, Michael R\"ockner, at the Departments of Mathematics of
Trento University, respectively, Bielefeld University, for their warm
hospitality during our stay in Trento, respectively, Bielefeld.

The second author 
gratefully acknowledges the warm hospitality of the Institute of
Applied Mathematics of the University of Bonn.


%

\printaddresses

\end{document}